\documentclass[10pt,reqno]{amsart}
\usepackage{amsmath,amsopn,amssymb,amsthm,multicol}
\usepackage{hyperref}
\usepackage{color}
\usepackage{graphicx}

\DeclareMathOperator{\diag}{diag}
\DeclareMathOperator{\Ad}{Ad}
\DeclareMathOperator{\ad}{ad}
\DeclareMathOperator{\Aut}{Aut}

\DeclareMathOperator{\tr}{tr}

\DeclareMathOperator{\Ric}{Ric}

\DeclareMathOperator{\Mark}{Mrk}

\newcommand{\fr}{\mathfrak}
\newcommand{\al}{\alpha}
\newcommand{\be}{\beta}

\newcommand{\bb}{\mathbb}
\DeclareMathOperator{\rnk}{rk}
\DeclareMathOperator{\SO}{SO}

\DeclareMathOperator{\Sp}{Sp}
\DeclareMathOperator{\SU}{SU}
\DeclareMathOperator{\U}{U}

\DeclareMathOperator{\F}{F}
\DeclareMathOperator{\E}{E}
\newcommand{\thickline}{\noalign{\hrule height 1pt}}

 \newtheorem{lemma} {Lemma} [section]
\newtheorem{theorem}[lemma]{Theorem} 
\newtheorem{remark}[lemma] {Remark} 
\newtheorem{prop} [lemma]{Proposition}  
\newtheorem{definition}[lemma] {Definition}

\newtheorem{conjecture}[lemma] {Conjecture}

\begin{document}

\title{Ricci flow on certain homogeneous spaces} 
\author{Marina Statha}
\address{University of Patras, Department of Mathematics, GR-26500 Rion, Greece}
\email{statha@math.upatras.gr} 
\medskip

\begin{abstract}
We study the behavior of the normalized Ricci flow of invariant Riemannian homogeneous metrics at infinity for  generalized Wallach spaces, generalized flag manifolds with four isotropy summands and second Betti number equal to one, and the Stiefel manifolds $V_2\bb{R}^n$ and $V_{1+k_2}\bb{R}^{n}$, with $n = 1+k_2+k_3$.  We use techniques from the theory of differential equations, in particular the Poincar\'e compactification.  


\end{abstract}

\maketitle
 

\section{Introduction}
\markboth{Marina Statha}{Ricci flow on certain homogeneous spaces}  

The Ricci flow equation was introduced by Hamilton in 1982 (\cite{Ha}), and is defined by
\begin{equation}\label{RicciFlow}
\frac{\partial g}{\partial t} = -2\Ric_g,
\end{equation}
where $g = g(t)$ is a curve on the space of Riemannian metrics $\mathcal{M}$ on a smooth manifold $M^n$ and $\Ric_g$ is the Ricci tensor of the Riemannian metric $g$.  The solution of this equation, the so called Ricci flow, is a 1-parameter family of metrics $g(t)$ in $M^n$. Intuitively, this is the heat equation for the metric $g$.

The Ricci flow (\ref{RicciFlow}) in general does not preserve the volume.  In the case of a compact manifold $M^n$ we consider the {\it normalized Ricci flow}
\begin{equation}\label{NormalRicci}
\frac{\partial g}{\partial t} = -2\Ric_g + \frac{2 r}{n}g,
\end{equation}
where $r = r(g(t)) = \int_{M}S_g du_g/ \int_{M} du_g$, $du_g$ is the volume element of $g$, and $S_g$ denotes the scalar curvature function of $g$.  Under this normalized flow, the volume of the solution metric is constant in time.  Equations (\ref{RicciFlow}) and (\ref{NormalRicci}) can be shown to be equivalent by reparametrizing time $t$ and by scaling the metric in space by a function of $t$.  The Einstein metrics, that is the Riemannian metrics of constant Ricci curvature which satisfy $\Ric_g = \lambda g$ (from now on call it Einstein equation) for some constant $\lambda\in\bb{R}$, are in this case the fixed points of the normalized Ricci flow (\ref{NormalRicci}).  In general, the Einstein equation reduces to a second order PDE and general existence results are difficult to be obtained.  Some methods are described in \cite{Bom}, \cite{BWZ} and \cite{WZ}.  Besides the detailed exposition on Einstein manifolds  in \cite{Be}, we refer to  \cite{W1}, \cite{W2} and \cite{Arv3} for more recent results.  For the case of homogeneous spaces $G/H$ the problem of finding all invariant Einstein metrics becomes slightly more accesible, due to the possibility of making symmetry assumptions, but  still it is not easy.  A special class of homogeneous spaces for which Einstein metrics have been completely classified are the generalized flag manifolds $G/K$ (of a compact simple Lie group $G$) with two (\cite{ArCh1}) three (\cite{Arv2}, \cite{Ki}) four and five isotropy summands (\cite{ArCh2}, \cite{ArChSa1}).  There is also some classification of Einstein metrics on flag manifolds with six isotropy summands (see in \cite{Arv3} for more details).  Also, the classification of invariant Einstein metrics for another class of homogeneous space, the generalized Wallach spaces, was only recently achieved (\cite{CN}).    
The problem becomes more difficult in case where the isotropy representation $\chi : H \to \Aut(T_{o}G/H)$, $(o=eH)$ of a homogeneous space $G/H$ contains equivalent summands.  This happens for example for Stiefel manifolds $V_{k}\bb{F}^{n}$, where $\bb{F}\in\{\bb{R}, \bb{C}, \bb{H}\}$.  In the papers \cite{ArSaSt1}, \cite{ArSaSt2} and \cite{ArSaSt3}, A. Arvanitoyergos, Y. Sakane and the author found (by using a technique which is described in detail in \cite{St}), Einstein metrics for several classes of Stiefel manifolds.   

An important property of the Ricci flow is that it preserves symmetries of the initial metric $g$.  This is due to the fact that the Ricci tensor is invariant under diffeomorphisms of the manifold $M$.  In general, the normalized Ricci flow (\ref{NormalRicci}) for an arbitrary manifold is a non-linear system of PDEs.  When restricted to the set of invariant metrics, such system reduces to an non-linear system of ODEs.  For this reason it is natural to study the Ricci flow on homogeneous spaces, that is a Riemmanian manifold $(M, g)$ with a closed subgroup $G$ of the isometries ${\rm Iso}(M, g)$, such that for any $p$ and $q$ in $M$, there exists
a $g\in G$ with $g(p) = q$.  In this case $M = G/H$, where $H = \{g\in G : g\cdot p = p\}$ is the isotropy subgroup at the point $p\in M$.  On such spaces we work with $G$-invariant metrics, i.e. metrics for which the map $\tau_\al : G/H \to G/H$, $gH \mapsto \al gH$ is an isometry.  It is natural to proceed the study of the Ricci flow using tools from the theory of dynamical systems.  
As mentioned above, to show existence of Einstein metrics is a not an easy task.  The use of the normalized Ricci flow on homogeneous spaces towards a qualitative study of homogeneous invariant Einstein metrics, has been used by various works, such as \cite{Ab}, \cite{AANS}, \cite{AbNi}, \cite{AnCh}, \cite{BW}, \cite{Bu}, \cite{GlPa},  \cite{GrMa1} and \cite{GrMa2}. 

In this paper we study the normalized Ricci flow of invariant metrics on certain homogeneous spaces with three and four isotropy summands, such as generalized Wallach spaces (that is a homogeneous space $G/H$ for which the tangent space $T_{o}(G/H)\cong\fr{m}$ written as direct sum of three summands $\fr{m}_i$, $i=1,2,3$ with the property $[\fr{m}_i, \fr{m}_i]\subset\fr{h}$), certain Stiefel manifolds $V_{k}\bb{R}^n$ (the set of all orthonormal $k$-frames on $\bb{R}^n$) and generalized flag manifolds (this is  homogeneous spaces $M = G/K$ where $G$ is a compact semisimple Lie group and $K$ is the cetralizer of a torus in $G$.  Equivalently, these are diffeomorphic to the adjoint orbit $\Ad(G)w$, for some $w\in\fr{g}$, the Lie algebra of $G$).  On such spaces the normalized Ricci flow (\ref{NormalRicci}) is equivalent to a homogeneous system of differential equations in $\bb{R}^3$ and $\bb{R}^4$.  So in order to study the behavior of such systems 
at infinity, we will use a method introduced by Poincar\'e, the so called {\it Poincar\'e compactification}.  This method allows us to study global phase portraits for polynomial systems.

The main contribution of the present work is that by using the Poincar\'e compactification we not only confirm previously obtained Einstein metrics as fixed points of dynamical systems deduced from normalized Ricci flow, but we also detect new homogeneous Einstein metrics on certain spaces.  The main theorems are the following:

\noindent{\bf Theorem A.}
{Let $G/H$ be a generalized Wallach space.  The normalized Ricci flow of a $G$-invariant Riemannian metric on $G/H$ has a finite number of singularities at infinity (in fact one to four). 
These fixed points determine the invariant Einstein metrics of $G/H$.}

\noindent
We also give some further analysis for the generalized Wallach space $\SO(k+l+m)/\SO(k)\times\SO(l)\times\SO(m)$, for $l=m$.

\noindent{\bf Theorem B.}
{(1) Let $V_2\bb{R}^n = \SO(n)/\SO(n-2)$ be a Stiefel manifold.  The normalized Ricci flow of the $\SO(n)$-invariant Riemannian metrics on $V_2\bb{R}^n$ has exactly one singularity at infinity.  This corresponds (up to scale) to the unique invariant Einstein metric.

\noindent(2) Let the Stiefel manifold $V_{1+k_2}\bb{R}^n \cong\SO(n)/\SO(k_3)$, with $n = 1+k_2+k_3$. The normalized Ricci flow on the space of invariant Riemannian metrics on $V_{1+k_2}\bb{R}^{1+2k_2}$, for $k_2 \geq 6$ has at least four singularities at infinity.  
These fixed points correspond (up to scale) to new $G$-invariant Einstein metrics on $V_{1+k_2}\bb{R}^{1+2k_2}$. 

\noindent(3) Let $G/H$ be the Stiefel manifold $V_5\bb{R}^7$.  The normalized Ricci flow on the space of invariant Riemannian metrics on $G/H$ possesses exactly four singularities at infinity.  These fixed points correspond (up to scale) to the $G$-invariant Einstein metrics on $G/H$. 
}

We are also motivated to state the following:

\noindent{\bf Conjecture 1.} 
{Let $G/H$ be the Stiefel manifold $V_{1+k_2}\bb{R}^{n}$, with $n = 1+k_2+k_3$.  Then, for $k_2 \geq 4$ and $k_3 > 1$, $G/H$ has four singularities at infinity, and for $k_2\geq 6$ and $k_3 = 2$  it has six singularities at infinity.  These fixed points correspond to the $G$-invariant Einstein metrics on $G/H$.}

\noindent{\bf Theorem C.}
{Let $G/K$ be a generalized flag manifold with four isotropy summands and $b_2(G/K) = 1$.  The normalized Ricci flow of $G$-invariant Riemannian metrics on $G/K$ has, for the case of exceptional flag manifold $\F_4$, $\E_7$ and $\E_8(\al_6)$, exactly three singularities at infinity.  For the case of flag manifold corresponding to $\E_8(\al_3)$ it has exactly five singularities at infinity.  These fixed points determine explicitly the three and five (up to scale) invariant Einstein metrics of $G/K$, respectively.  
}



\section{The normalized Ricci flow}

\subsection{Ricci tensor and scalar curvature}

Let $G$ be a compact semisimple Lie group, $H$ a connected closed subgroup of $G$ and let  $\frak g$ and $\fr{h}$  be the corresponding Lie algebras.  The Killing form $B$ of $\frak g$ is negative definite, so we can define an $\Ad(G)$-invariant inner product $-B$ on $\frak g$.  Let $\frak g$ = $\frak h \oplus
\frak m$ be a reductive decomposition of $\frak g$ with respect to $-B$ so that $\left[\,\frak h,\, \frak m\,\right] \subset \frak m$ and $\frak m\cong T_o(G/H)$ where $o$ is the identity coset of $G/H$.  Any $G$-invariant metric $g$ on $G/H$ corresponds to an $\Ad(H)$-invariant inner product $\langle\cdot, \cdot\rangle$ on $\fr{m}$ and vice versa.  Let $\{X_j\}$ be a $\langle\cdot, \cdot\rangle$-orthonormal basis of $\frak m$.  The Ricci tensor  $\Ric_g$ of the metric $g$ is given as follows (\cite[p. 185]{Be}):
\begin{eqnarray}\label{r}
&&\Ric_g(X, Y)=-\frac12\sum_i\langle [X, X_i], [Y, X_i]\rangle +\frac12B(X, Y)
       +\frac14 \sum _{i, j}\langle [X_i, X_j], X\rangle\langle [X_i, X_j], Y\rangle.\label{ric1}
\end{eqnarray} 
The scalar curvature $S_g = \tr\Ric_g$ of $g$ is given by (\cite[p. 186]{Be}):
\begin{eqnarray}
&&S_g = \frac{1}{4}\sum_{i,j} |[X_i, X_j]_{\fr{m}}|^{2} -\frac{1}{2}\sum_i B(X_i, X_i).\label{scalar}
\end{eqnarray}

If the isotropy representation of $G/H$ is decomposed into a sum of mutually non equivalent irreducible summands,
then we will also use the following alternative expression for the Ricci tensor. Let
 $
{\frak m} = {\frak m}_1 \oplus \cdots \oplus {\frak m}_q
$
be a decomposition into mutually non equivalent irreducible $\Ad(H)$-modules.  
Then any $G$-invariant metric on $G/H$ can be expressed as follows:  
\begin{eqnarray}
 \langle\cdot, \cdot\rangle  = x_1(-B)|_{\fr{m}_1} + x_2(-B)|_{\fr{m}_2} + \cdots + x_{q}(-B)|_{\fr{m}_q},  \label{eq2}
\end{eqnarray}
for positive real numbers $(x_1, \dots, x_q)\in\bb{R}^{q}_{+}$.  Note that  $G$-invariant symmetric covariant 2-tensors on $G/H$ are of the same form as the Riemannian metrics (although they  are not necessarilly  positive definite).  In particular, the Ricci tensor $\Ric_g$ of a $G$-invariant Riemannian metric on $G/H$ is of the same form as (\ref{eq2}), that is 
$
\Ric_g = y_1(-B)|_{\fr{m}_1} +  y_2(-B)|_{\fr{m}_2} + \cdots +  y_q(-B)|_{\fr{m}_q}
$
for some $y_i\in\bb{R}, i=1,2,\ldots,q$.  

Let $\lbrace e_{\alpha}^{(k)} \rbrace_{\alpha=1}^{d_k}$, where $d_k = \dim\fr{m}_k$, be a $(-B)$-orthonormal basis of $\fr{m}_k$.  Then the set $\lbrace X_{\alpha}^{(k)} = e_{\alpha}^{(k)}/\sqrt{x_{k}}\rbrace$ is a $\langle\cdot,\cdot\rangle$-orthonormal basis of $\fr{m}_k$.  If we denote by $r_k = \Ric_{g}(X_{\alpha}^{(k)}, X_{\alpha}^{(k)})$, then we obtain $r_k = ({1}/{x_k})$ $\Ric_{g}(e_{\alpha}^{(k)},$ $e_{\alpha}^{(k)})$ that is $\Ric_{g}(e_{\alpha}^{(k)}, e_{\alpha}^{(k)}) = x_{k}r_{k}$.  Thus the Ricci tensor is written as 
$
\Ric_g = x_1r_1(-B)|_{\fr{m}_1} +  x_2r_2(-B)|_{\fr{m}_2} + \cdots +  x_qr_q(-B)|_{\fr{m}_q},
$
where the $r_i$'s are the components of the Ricci tensor on each $\fr{m}_i$ for $i = 1,2,\ldots,k$. 
Now let ${A^\gamma_{\alpha \beta}}= -B ([e_{\alpha}^{(i)}, e_{\beta}^{(j)}], e_{\gamma}^{(k)})$ so that
$[e_{\alpha}^{(i)}, e_{\beta}^{(j)}]
= \displaystyle{\sum_{\gamma}
A^\gamma_{\alpha \beta} e_{\gamma}^{(k)}}$ and set $A_{ijk}: = \displaystyle{k \brack {ij}}=\sum (A^\gamma_{\alpha \beta})^2$, where the sum is taken over all indices $\alpha, \beta, \gamma$ with $e_\alpha^{(i)} \in {\frak m}_i,\ e_\beta^{(j)} \in {\frak m}_j,\ e_\gamma^{(k)} \in {\frak m}_k$ (cf. \cite{WZ}).  Then the positive numbers $A_{ijk}$ are independent of the 
$B$-orthonormal bases chosen for ${\frak m}_i, {\frak m}_j, {\frak m}_k$, and 
$A_{ijk} = A_{jik} = A_{kij}$.  
We have the following:

\begin{lemma}\textnormal{(\cite{PaSa})}\label{ParkSakane}
Let $G/H$ be a homogeneous space where $G$ is compact and semisimple Lie group.  Let $g$ the $G$-nvariant metric on $G/H$ given by the $\Ad(H)$-invariant inner products (\ref{eq2}).  Then:
\begin{itemize}
\item[(1)] The components $r_1,\ldots,r_{s}$ of the Ricci tensor $\Ric_g$ are given by 
\begin{equation}
{r}_k = \frac{1}{2x_k}+\frac{1}{4d_k}\sum_{j,i}
\frac{x_k}{x_j x_i} A_{jik}
-\frac{1}{2d_k}\sum_{j,i}\frac{x_j}{x_k x_i} A_{kij}
 \quad (k= 1,\ \dots, q)    
\end{equation}
where the sum is taken over $i, j =1,\dots, q$.
\item[(2)] The scalar curvature $S_g = {\rm tr}\Ric_g = \sum_{i=1}^{s}d_ir_i$ is given by 
\begin{equation}
S_g = \frac{1}{2}\sum_{i=1}^{s}\frac{d_i}{x_i} - \frac{1}{4}\sum_{i,j,k}\frac{x_k}{x_ix_j}A_{ijk} 
\end{equation}
\end{itemize}
\end{lemma}

\begin{remark}
\textnormal{Since by assumption the tangent space $\fr{m}$ of $G/H$ decomposes into $\Ad(H)$-modules $\fr{m}_{i}, \fr{m}_{j}$ which are mutually non equivalent for any $i\neq j$, it is $\Ric_{g}(\fr{m}_{i}, \fr{m}_{j})=0$ whenever $i\neq j$.  Thus, by Lemma \ref{ParkSakane} it follows that $G$-invariant Einstein metrics on $G/H$ are exactly the positive real solutions $g=(x_1, \ldots, x_q)\in\bb{R}^{q}_{+}$  of the  polynomial system $\{r_1=\lambda, r_2=\lambda, \ldots, r_{q}=\lambda\}$, where $\lambda\in \bb{R}_{+}$ is the Einstein constant.
}
\end{remark}

\subsection{The normalized Ricci flow}
Let $(M = G/H, g)$ a Riemannian homogeneous space and let $\mathcal{M}_{1}^{G}$ be the set of $G$-invariant metrics with total volume 1, that is $\int_{M}dv_{g} = 1$.  Then for every $g$ in $\mathcal{M}_{1}^{G}$ its scalar curvature $S_g$ is a constant function on $M$.  Therefore we have $r = \int_{M}S_g dv_{g}/\int_{M}dv_g = S_g$.  Thus, for a $G$-invariant metric on $G/H$ the normalized Ricci flow (\ref{NormalRicci}) is equivalent to 
$$
\frac{\partial g}{\partial t} = -2\Ric_g + \frac{2S_g}{n}g,
$$
where $n$ is the dimension of $G/H$.  Actually for the $G$-invariant metric (\ref{eq2}) the normalized Ricci flow reduces to the following system of ODEs:
\begin{equation}\label{normalRicci}
\Big\{ \dot{x}_1 =2x_1 r_1 + {\displaystyle\frac{2 S_g}{n}x_1},\ \ \dot{x}_2 =2x_2 r_2 + {\displaystyle\frac{2S_g}{n}x_2},\ \ldots, \ \dot{x}_q =2x_q r_q + {\displaystyle\frac{2 S_g}{n}x_q}\Big\}.
\end{equation}

\section{Generalized Wallach spaces}

Let $G/H$ be a reductive homogeneous space with $G$ a compact and semisimple Lie group and $H$ a compact subgroup of $G$.  Let $\fr{g} = \fr{h}\oplus\fr{m}$ the reductive decomposition of $G/H$, that is $\Ad(H)\fr{m}\subset \fr{m}$, with $\fr{m}\cong T_{o}(G/H)$.  Then $G/H$ is called a {\it generalized Wallach space} if the module $\fr{m}$ decomposes into a direct sum of three $\Ad(H)$-invariant irreducible modules pairwise orthogonal with respect to $-B$ (the Killing form of $G$), i.e. $\fr{m} = \fr{m}_1\oplus\fr{m}_2\oplus\fr{m}_3$, such that $[\fr{m}_i, \fr{m}_i]\subset\fr{h}$, for $i = 1,2,3$.  Every generalized Wallach space admits a three parameter family of $G$-invariant Riemannian metrics determined by $\Ad(H)$-invariant inner products:
\begin{equation}\label{Nikonorov}
\langle\cdot, \cdot\rangle = x_1(-B)|_{\fr{m}_1} + x_2(-B)|_{\fr{m}_2} + x_3(-B)|_{\fr{m}_3}, \  x_i\in\bb{R}_{+}, i =1,2,3,
\end{equation}  
where $x_1, x_2, x_3$ are positive real numbers.  We will denote such metrics with $g = (x_1, x_2, x_3)$.

The classification of generalized Wallach spaces $G/H$ was obtained in \cite{Ni1} and \cite{ChKaLi}:

\begin{theorem}
Let $G/H$ be a connected and simply connected compact homogeneous space.  Then $G/H$ is a generalized Wallach space if and only if it is one of the following types:

\smallskip
\noindent(1) $G/H$ is a direct product of three irreducible symmetric spaces of compact type.

\smallskip
\noindent(2) The group $G$ is simple and the pair $(\fr{g}, \fr{k})$ is one of the pairs in Table 1.

\smallskip
\noindent(3) $G = F \times F \times F \times F$ and $H = \diag(F) \subset G$ for some connected, compact, simple Lie group $F$, with the following description on the Lie algebra level: 
$(\fr{g}, \fr{h}) = (\fr{f}\oplus\fr{f}\oplus\fr{f}\oplus\fr{f}, \diag(\fr{f}) = \{(X,X,X,X) | X\in\fr{f}\})$, where $\fr{f}$ is the Lie algebra of $F$, and (up to permutation) 
$\fr{m}_1 = \{(X,X,-X,-X) | X\in\fr{f}\}$, $\fr{m}_2 = \{(X,-X,X,-X) | X\in\fr{f}\}$, $\fr{m}_3 = \{(X,-X,-X,X) | X\in\fr{f}\}$.
\end{theorem} 

\smallskip
 \begin{center}
{\small {\bf Table 1.} The pairs $(\fr{g}, \fr{h})$ corresponding to generalized Wallach spaces $G/H$ with $G$ simple }  
\end{center} 
{\small\small \begin{center}
 \begin{tabular}{l|l|l|l|l|l}
                \thickline
\mbox{\bf{GWS.}} & $\fr{g}$         &     $\fr{h}$     &   $d_1$   &   $d_2$   &   $d_3$        \\
\hline
$\bf{1}$ & $\fr{so}(k+l+m)$  &  $\fr{so}(k)\oplus\fr{so}(l)\oplus\fr{so}(m)$  &  $kl$   &  $km$   &  $lm$\\
$\bf{2}$ & $\fr{su}(k+l+m)$  &  $\fr{su}(k)\oplus\fr{su}(l)\oplus\fr{su}(m)$  &  $2kl$  &  $2km$  &  $2lm$ \\
$\bf{3}$ & $\fr{sp}(k+l+m)$  &  $\fr{sp}(k)\oplus\fr{sp}(l)\oplus\fr{sp}(m)$  &  $4kl$  &  $4km$  &  $4lm$\\
$\bf{4}$ & $\fr{su}(2l),\, l\geq 2$  & $\fr{u}(l)$                       &  $l(l-1)$   & $l(l+1)$ &  $l^2-1$  \\
$\bf{5}$ & $\fr{so}(2l), \, l\geq 4$ & $\fr{u}(1)\oplus\fr{u}(l-1)$      &  $2(l-1)$   &  $2(l+1)$ &  $(l-1)(l-2)$\\
$\bf{6}$ & $\fr{e}_6$       &   $\fr{su}(4)\oplus2\fr{sp}(1)\oplus\bb{R}$    &  $16$     &  $16$   &  $24$\\
$\bf{7}$ & $\fr{e}_6$       &   $\fr{so}(8)\oplus\bb{R}^2$                   &  $16$     &  $16$   &  $16$\\
$\bf{8}$ & $\fr{e}_6$       &   $\fr{sp}(3)\oplus\fr{sp}(1)$                 &  $14$     &  $28$   &  $12$\\
$\bf{9}$ & $\fr{e}_7$       &   $\fr{so}(8)\oplus3\fr{sp}(1)$                &  $32$     &  $32$   &  $32$\\
$\bf{10}$ & $\fr{e}_7$       &   $\fr{su}(6)\oplus\fr{sp}(1)\oplus\bb{R}$     &  $30$     &  $40$   &  $24$\\
$\bf{11}$ & $\fr{e}_7$       &   $\fr{so}(8)$                                 &  $35$     &  $35$   &  $35$\\
$\bf{12}$ & $\fr{e}_8$       &   $\fr{so}(12)\oplus2\fr{sp}(1)$               &  $64$     &  $64$   &  $48$\\  
$\bf{13}$ & $\fr{e}_8$       &   $\fr{so}(8)\oplus\fr{so}(8)$                 &  $64$     &  $64$   &  $64$\\ 
$\bf{14}$ & $\fr{f}_4$       &   $\fr{so}(5)\oplus2\fr{sp}(1)$                &  $8$      &  $8$    &  $20$\\  
$\bf{15}$ & $\fr{f}_4$       &   $\fr{so}(8)$                                 &  $8$      &  $8$    &  $8$     
\\
  \thickline
 \end{tabular}
 \end{center} }
\medskip
\smallskip

\subsection{Ricci tensor for generalized Wallach spaces}

Let $d_i = \dim\fr{m}_i$, $i = 1,2,3$.  From the definition of $A_{ijk}$ it is easy to see that for the generalized Wallach spaces, the $A_{ijk} = 0$ if two of the indicies are equal.  Therefore, we only need to compute the number $A_{123}$.  It is $d_i \geq 2A_{123}$ for any $i= 1,2,3$ (see \cite{Ni2}).  Hence the constants $a_i = A_{123}/d_i$, $i\in\{1,2,3\}$ are such that $(a_1, a_2, a_3)\in (0, 1/2]^3$.  In Table 2 we give these numbers for the generalized Wallach spaces listed in Table 1 (these numbers were computed in \cite{Ni1} and \cite{ChKaLi}) 

\smallskip
 \begin{center}
{\small {\bf Table 2.} The numbers $a_i$, $i=1,2,3$ for the generalized Wallach spaces $G/H$ with $G$ simple }  
\end{center} 
{\small\small \begin{center}
 \begin{tabular}{l|l|l|l||l|l|l|l}
                \thickline
   \mbox{\bf{GWS.}}      &   $a_1$   &   $a_2$   &   $a_3$ &  \mbox{\bf{GWS.}} &  $a_1$ & $a_2$   &   $a_3$     \\
\hline
  $\bf{1}$  &  $m/2(k+l+m-2)$   &  $l/2(k+l+m-2)$   &  $k/2(k+l+m-2)$ &  $\bf{9}$  &  $2/9$            &  $2/9$            &  $2/9$\\
  $\bf{2}$  &  $m/2(k+l+m)$     &  $l/2(k+l+m)$     &  $k/2(k+l+m)$ &  $\bf{10}$ &  $2/9$            &  $1/6$            &  $5/18$\\
  $\bf{3}$  &  $m/2(k+l+m+1)$   &  $l/2(k+l+m+1)$   &  $k/2(k+l+m+1)$ &  $\bf{11}$ &  $5/18$           &  $5/18$           &  $5/18$\\
  $\bf{4}$  &  $(l+1)/4l$       &  $(l-1)/4l$       &  $1/4$  &  $\bf{12}$ &  $1/5$            &  $1/5$            &  $4/15$\\  
  $\bf{5}$  &  $(l-2)/4(l-1)$   &  $(l-2)/4(l-1)$   &  $1/2(l-1)$ & $\bf{13}$ &  $4/15$           &  $4/15$           &  $4/15$\\ 
  $\bf{6}$  &  $1/4$            &  $1/4$            &  $1/6$ & $\bf{14}$ &  $5/18$           &  $5/18$           &  $1/9$\\  
  $\bf{7}$  &  $1/6$            &  $1/6$            &  $1/6$ &  $\bf{15}$ &  $1/9$            &  $1/9$            &  $1/9$  
\\  $\bf{8}$  &  $1/4$            &  $1/8$            &  $7/24$ & & & &  
\\
  \thickline
 \end{tabular}
 \end{center} }
\medskip
\smallskip

From Lemma \ref{ParkSakane} the components of the Ricci tensor for the metric which corresponds to the $\Ad(H)$-invariant inner products (\ref{Nikonorov}), are given as follows
$$
r_i = \frac{1}{2 x_i} + \frac{a_{i}}{2}\left(\frac{x_i}{x_j x_k} - \frac{x_k}{x_i x_j} - \frac{x_j}{x_i x_k}\right),
$$  
where $i, j, k\in\{1,2,3\}$ with $i\neq j\neq k\neq i$.
The scalar curvature is given by the following 
$$
S_{g} = \frac{1}{2}\left(\frac{d_1}{x_1} + \frac{d_2}{x_2} + \frac{d_3}{x_3}\right) - A_{123}\left(\frac{x_1}{x_2x_3} + \frac{x_2}{x_1x_3} + \frac{x_3}{x_1x_2}\right).
$$

The $G$-invariant metric on a generalized Wallach space $G/H$ corresponding to the inner product (\ref{Nikonorov}) is Einstein if and only if $r_1 = r_2 = r_3$.  This is equivalent to the polynomial system
\begin{eqnarray}
(a_2 + a_3)(a_1x_2^2 + a_1x_3^2 - x_2x_3) + (a_2x_2 + a_3x_3) - (a_1a_2 + a_1a_3 + 2a_2a_3)x_1^2 = 0 \nonumber\\
(a_1 + a_3)(a_2x_1^2 + a_2x_3^2 - x_1x_3) + (a_1x_1 + a_3x_3) - (a_1a_2 +2a_1a_3 + a_2a_3)x_2^2 = 0.
\end{eqnarray} 

The first general result about the number of invariant Einstein metrics on generalized Wallach spaces is the following:

\begin{theorem}{\textnormal{(\cite{LoNiFi})}}\label{general}
Let $G/H$ be a generalized Wallach space  with pairwise non equivalent submodules $\fr{m}_i$.  Then $G/H$ admits from one to four invariant Einstein metrics (up to scale and homothety).
\end{theorem}

We take the $G$-invariant metric of the form $g = (1, y_1, y_2)$, where $y_i = x_i/x_1$, $i = 1,2$ on generalized Wallach spaces.  Then from  \cite{LoNiFi} and \cite{ChKaLi} we have the following Einstein metrics for  {\small\small\bf GWS} corresponding to exceptional groups.

\smallskip
 \begin{center}
{\small {\bf Table 3.} The Einstein metrics for the generalized Wallach spaces $G/H$ }  
\end{center} 
{\small\small \begin{center}
 \begin{tabular}{l|l|l|l|l}
                \thickline
 \mbox{\bf{GWS.}}        & $g_1 = (1, y_2, y_3)$ & $g_2 = (1, y_2, y_3)$ & $g_3 = (1, y_2, y_3)$ & $g_4 = (1, y_2, y_3)$        \\
         \hline
$\bf{6}$  &  $(1, 0.6, 0.8)$   &  $(1, 1.66667, 1.33333)$   &  $-$ & $-$  \\
\hline
$\bf{7}$  & $(1, 1, 1)$  &  $(1, 0.5, 0.5)$  &  $(1, 2, 1)$  &  $(1, 1, 2)$\\
\hline
$\bf{8}$ &   $(1, 1.4618, 1.88845)$ & $(1, 0.8640, 0.4838)$  & $-$  & $-$\\
\hline
$\bf{9}$ &   $(1, 1, 1)$ &  $(1, 1.25, 1)$  & $(1, 0.8, 0.8)$  & $(1, 1, 1.25)$\\
\hline
$\bf{11}$  &  $(1, 1, 1)$ & $(1, 1, 0.8)$  & $(1, 1.25, 1.25)$  & $(1, 0.8, 1)$\\
\hline
$\bf{12}$  &  $(1, 1, 1.45608)$  &  $(1, 1, 0.68677)$  &  $-$  & $-$\\
\hline
$\bf{13}$  & $(1, 1, 1)$  & $(1, 1, 0.875)$  & $(1, 1.14285, 1.14285)$ & $(1, 0.875, 1)$\\
\hline
$\bf{14}$  &  $(1, 0.4852, 0.8251)$  &  $(1, 2.0606, 1.700349)$  & $-$  & $-$\\
\hline
$\bf{15}$  & $(1, 1, 1)$ & $(1, 3.5, 1)$ & $(1, 1, 3.5)$  & $(1, 0.28571, 0.28571)$
\\
  \thickline
 \end{tabular}
 \end{center} }
\medskip
\smallskip


The above list does not include case  {\small\small\bf GWS.10}, that has been studied in \cite{AnCh}.  
Note that the spaces  {\small\small\bf GWS.2} and {\small\small\bf GWS.5} are also generalized flag manifolds with three isotropy summands (this is also true for spaces {\small\small\bf GWS.7} and {\small\small\bf GWS.10}).  The classification of Einstein metrics on these flag manifolds was given in \cite{Ki} and \cite{Arv2}, which agrees with the following theorem:
\begin{theorem}{\textnormal{(\cite{LoNiFi})}}
If a generalized Wallach space $G/H$ with pairwise non-isomorphic modules $\fr{m}_i$ satisfy the equality $a_1+a_2+a_3 = 1/2$, then $G/H$ admits four families of proportional invariant Einstein metrics.  These metrics have the form
\begin{eqnarray}
(1)\ ((1-2a_1)q, (1-2a_2)q, 2(a_1+a_2)q), && (2)\ ((1-2a_1)q, (1-2a_2)q, 2(1-a_1-a_2)q), \nonumber\\
(3)\ ((1-2a_1)q, (1+2a_2)q, 2(a_1+a_2)q), && (4)\ ((1+2a_1)q, (1-2a_2)q, 2(a_1+a_2)q), 
\end{eqnarray}
where $q\in\bb{R}$.
\end{theorem}

The following theorem gives a good insight about the number of invariant Einstein metrics on the generalized Wallach space $\SO(k+l+m)/\SO(k)\times\SO(l)\times\SO(m)$.
\begin{theorem}{\textnormal{(\cite{CN})}}\label{china}
Assume that $1\le k\le l\le m$ and $\l\ge 2$ and let $G/H=\SO(k+l+m)/\SO(k)\times\SO(l)\times\SO(m)$. Then,
for $k>\sqrt{2m+2l-4}$ the number of invariant Einstein metrics on $G/H$ is $4$, and for $k<\sqrt{m+l}$  the number of invariant Einstein metrics on $G/H$ is $2$.
\end{theorem}

Finally, for the generalized Wallach spaces  {\small\small\bf GWS.3} and {\small\small\bf GWS.4}, the Einstein metrics are given as solutions of equation (9) in the paper \cite[p. 51]{LoNiFi}.  
We give for {\small\small\bf GWS.3} some examples of Einstein metrics (up to scale):




\smallskip
\noindent{\small\small\bf GWS.3a} For $\Sp(6)/\Sp(1)\times\Sp(2)\times\Sp(3)$ we have the following metrics $g_i = (1, x_1, x_2), i=1,2,3,4$
$
\centerline{ (1, 0.38050, 0.46780),  (1, 1.23251, 1.39606),
(1, 3.26361, 1.60389),  (1, 1.30670, 3.18223).}
$

\smallskip
\noindent{\small\small\bf GWS.3b} For $\Sp(14)/\Sp(2)\times\Sp(5)\times\Sp(7)$ we have the following metrics $g_i = (1, x_1, x_2), i=1,2,3,4$
$
\centerline{(1, 0.40168, 0.52944), (1, 1.24716, 1.53155), 
(1, 2.94748, 1.67504),  (1, 1.27217, 2.71689).}
$




\section{Stiefel manifolds} 

We embed the group $\SO(n-k)$ in $\SO(n)$ as $\begin{pmatrix}
1_k & 0\\
0 & C
\end{pmatrix}$ where $C\in\SO(n-k)$.  The Killing form of $\fr{so}(n)$ is $B(X, Y)=(n-2)\tr XY$.
Then with respect to $-B$ the subspace $\fr{m} =\fr{so}(n-k)^{\perp}$ in $\fr{so}(n)$, may be identified with the set 
of matrices of the form
$
\left\lbrace \begin{pmatrix}
D_k & A\\
-A^t & 0_{n-k}
\end{pmatrix} : D_k \in \fr{so}(k), A\in M_{k\times(n-k)}(\bb{R}) \right\rbrace.
$
Let $E_{ab}$ denote the $n\times n$ matrix with $1$ at the $(ab)$-entry and $0$ elsewhere.  Then the set $\mathcal{B}=\{e_{ab}=E_{ab}-E_{ba}: 1\le a\le k,\ 1\le a<b\le n\}$
constitutes a $-B$-orthogonal basis of $\fr{m}$.  Note that $e_{ba}=-e_{ab}$, thus we have: 

\begin{lemma}\label{brac}
If all four indices are distinct, then the Lie brackets in $\mathcal{B}$ are zero.
Otherwise,
$[e_{ab}, e_{bc}]=e_{ac}$, where $a,b,c$ are distinct.
\end{lemma}

Next, we study the isotropy representation of $V_k\mathbb{R}^n=G/H=\SO(n)/\SO(n-k)$.  Let $\lambda _n$ denote the standard representation of $\SO(n)$ (given by the natural action of $\SO(n)$ on $\mathbb{R}^n$).  If $\wedge ^2\lambda _n$ denotes the second exterior power of $\lambda _n$, then $\Ad ^{\SO(n)}=\wedge ^2\lambda _n$.
The isotropy representation $\chi : \SO(n)\to \Aut(\fr{m})$ ($\fr{m}\cong T_{o}(G/H)$) of $G/H$ is characterized by the property
$\left.\Ad ^{\SO(n)}\right |_{\SO(n-k)}=\Ad ^{\SO(n-k)}\oplus\chi$. 
We compute
\begin{equation}\label{isotropy1}
\left.\Ad ^{\SO(n)}\right |_{\SO(n-k)}=\wedge ^2\lambda _n\big| _{\SO(n-k)}=\wedge ^2 (\lambda _{n-k}\oplus k)=\wedge ^2\lambda _{n-k}\oplus \wedge ^2k\oplus (\lambda _{n-k}\otimes k),
\end{equation}
where $k$ denotes the trivial $k$-dimensional representation.  Therefore, the isotropy representation is given by $\chi = 1\oplus \cdots\oplus 1\oplus\lambda _{n-k}\oplus \cdots\oplus\lambda _{n-k}$.  This decomposition induces an $\Ad(H)$-invariant decomposition of $\fr{m}$ given by $\fr{m}=\fr{m}_1\oplus\cdots\oplus\fr{m}_s$,
where the first ${k}\choose{2}$ $\Ad(H)$-modules are $1$-dimensional and the rest $k$ are $(n-k)$-dimensional.  It is clear that the isotropy representation of $V_k\bb{R}^n$ contains equivalent summands, so a complete description of all $G$-invariant metrics is rather hard.  In \cite{ADN1} the authors introduced a method for proving existence of homogeneous Einstein metrics by assuming additional symmetries.  In \cite{St} is presented a systematic and organized description of such metrics.

\subsection{The Stiefel manifold $V_2\bb{R}^{n}\cong \SO(n)/\SO(n-2)$}
The isotropy representation of $V_2\bb{R}^n$, is expressed as a direct sum $\chi = 1\oplus\chi_1\oplus\chi_2$, where $\chi_1\approx \chi_2= \lambda_{n-2}$ is the standard representation of $\SO(n-2)$.  This decomposition induces an $\Ad(\SO(n-2))$-invariant decomposition of $\fr{m}$ given by $\fr{m} = \fr{m}_{0}\oplus\fr{m}_1\oplus\fr{m}_2$.  
Even though an $\SO(n)$-invariant metric on $V_2\bb{R}^n$ depends on four parameters, it can be shown (cf. \cite{Ke}) that it can be descrided by an $\Ad(\SO(n-2))$-invariant inner product of $\fr{m}$ of the form:
\begin{equation}\label{metricSt}
\langle\cdot, \cdot\rangle = x_0(-B)|_{\fr{m}_0} + x_1(-B)|_{\fr{m}_1} + x_2(-B)|_{\fr{m}_2}, \  x_i\in\bb{R}_{+}, i =1,2,3.
\end{equation}
Therefore,  for the Ricci tensor of metrics corresponding to inner products (\ref{metricSt}) we can use the Lemma \ref{ParkSakane}.  By using Lemma \ref{brac} the only non-zero number is $A_{012}$ and equals to $1/2$.  Hence we have,

\begin{prop}
(1) The components of the Ricci tensor for the metric (\ref{metricSt}) are given as follows
\begin{equation}\label{ricciV2}
\left. 
{\small \begin{array}{ll}
r_0  =  \displaystyle{\frac{1}{2 {x_0}} - \frac{1}{4} \left(\frac{{x_1}}{{x_0}
   {x_2}}+\frac{{x_2}}{{x_0} {x_1}} - \frac{{x_0}}{{x_1} {x_2}} \right), } \ \  
r_1 = \displaystyle{\frac{1}{2x_1}  - \frac{1}{4(n-2)}\left(\frac{x_0}{x_1 x_2} + \frac{x_2}{x_0x_1} - \frac{x_1}{x_0 x_2}\right) } \\ \\
r_2 = \displaystyle{\frac{1}{2x_2}  -\frac{1}{4(n-2)}\left(\frac{x_0}{x_1 x_2} + \frac{x_1}{x_0 x_2} - \frac{x_2}{x_0 x_1}\right) } 
\end{array} } \right\}
\end{equation} 
(2) The scalar curvature $S_g$ is given by
{ \begin{eqnarray}\label{scalV2}
S_g &=& \frac{1}{2x_0} + \frac{1}{4}\left(-\frac{x_1}{x_0x_2}-\frac{x_2}{x_0x_1}-\frac{x_0}{x_1x_2}\right) + \frac{n-2}{2}\left(\frac{1}{x_2} + \frac{1}{x_1}\right)  
\end{eqnarray} }
\end{prop}

\begin{theorem}\textnormal{(\cite{Arv1}, \cite{Ke})}
The Stiefel manifold $V_2\bb{R}^n = \SO(n)/\SO(n-2)$ admits (up to scale) exactly one $\SO(n)$-invariant Einstein metric which is given explicitly as $(1, ({n-1})/{2(n-2)}, ({n-1})/{2(n-2)})$.
\end{theorem}

\subsection{The Stiefel manifolds $V_{1+k_2}\bb{R}^n$}
Let $G/H$ be the Stiefel manifold $V_{1+k_2}\bb{R}^n = \SO(n)/\SO(k_3)$, with $n = 1+k_2+k_3$.  The isotropy representation for this case according to (\ref{isotropy1}) contains equivalent summands.  Next, we describe a special class of invariant metrics on this space (for more details see for example \cite{St}, \cite{ArSaSt1} and \cite{ArSaSt2}).  The basic approach is to use an appropriate subgroup $K$ of $G$, such that the special class of $\Ad(K)$-invariant inner products, which are a subset of $\Ad(H)$-invariant inner products, are diagonal.  In order to have this, it is sufficient for the subgroup $K$ to satisfy the condition $H\subset K \subset N_{G}(H)\subset G$.

We take the subgroup $K = \SO(k_2)\times\SO(k_3)$ of $\SO(n)$.  Then, for the tangent space $\fr{m} \cong T_{o}(G/H)$, we consider the irreducible, $\Ad(K)$-invariant and non-equivalent decomposition: 
$
\fr{m} = \fr{so}(k_2)\oplus  \fr{m}_{12} \oplus  \fr{m}_{13} \oplus  \fr{m}_{23}\footnote{The direct sum $\fr{m}_{12} \oplus  \fr{m}_{13} \oplus  \fr{m}_{23}$ is the tangent space of generalized Wallach space { {\bf GWS.1}}}. 
$

Then the $G$-invariant metrics on $G/H$ determined by the  $\Ad(\SO(k_2)\times\SO(k_3))$-invariant scalar products on $\fr{m}$ are given by 
\begin{equation} \label{metric2} 
 \langle \cdot, \cdot \rangle = x_2 (-B) |_{ \fr{so}(k_2)} + x_{12} (-B) |_{ \fr{m}_{12}} + x_{13} (-B) |_{ \fr{m}_{13}} + x_{23} (-B) |_{ \fr{m}_{23}}
\end{equation}

Then by using Lemma \ref{brac} it follows that the only non zero triplets (up to permutation of indices) are 
$
A_{222},\ A_{2(12)(12)},\ A_{2(23)(23)},\ A_{(13)(12)(23)}, 
$
where $A_{iii}$ is non zero only for $k_2 \geq 3$.

\begin{lemma}\textnormal{(\cite{ArSaSt1})}\label{lemma5.20} For $a, b, c = 1, 2, 3$ and $(a - b)(b - c) (c - a) \neq 0$ the following relations hold:
\begin{equation*}\label{eq14}
\begin{array}{lll} 
 A_{aaa} = \displaystyle{\frac{k_a (k_a -1)(k_a -2)}{2 (n -2)}} ,   &  A_{a (a b) (a b)} = \displaystyle{\frac{k_a  k_b (k_a -1)}{2 (n -2)}}, &  A_{(a c) (a b ) (b c)} = \displaystyle{\frac{k_a  k_b  k_c}{2 (n -2)}}. 
\end{array} 
\end{equation*}
\end{lemma} 

\begin{lemma}\label{ricV4Rn}
(1) The components of the Ricci tensor $\Ric$ for the invariant metric $\langle\cdot, \cdot\rangle $ on $G/H$ defined by  {\em (\ref{metric2})}, are given as follows:  
\begin{equation}\label{eq19}
\left. {\small \begin{array}{l} 
 r_2  = 
\displaystyle{\frac{k_2-2}{4 (n -2) x_2} +
\frac{1}{4 (n -2)} \biggl(  \frac{x_2}{{x_{12}}^2} +k_3 \frac{x_2}{{x_{23}}^2} \biggr)},
\\  \\
r_{12}   =  \displaystyle{\frac{1}{2 x_{12}} +\frac{k_3}{4 (n -2)}\biggl(\frac{x_{12}}{x_{13} x_{23}} - \frac{x_{13}}{x_{12} x_{23}} - \frac{x_{23}}{x_{12} x_{13}}\biggr) -
\frac{1}{4 (n -2)} \biggl(  (k_2-1) \frac{x_2}{{x_{12}}^2} \biggr)},
\\  \\
r_{23}   =  \displaystyle{\frac{1}{2 x_{23}} +\frac{1}{4 (n -2)}\biggl(\frac{x_{23}}{x_{13} x_{12}} - \frac{x_{13}}{x_{12} x_{23}} - \frac{x_{12}}{x_{23} x_{13}}\biggr) -
\frac{1}{4 (n -2)} \biggl(  (k_2-1) \frac{x_2}{{x_{23}}^2} \biggr)}, 
\\  \\
r_{13}   =  \displaystyle{\frac{1}{ 2 x_{13}} +\frac{k_2}{4 (n -2)}\biggl(\frac{x_{13}}{x_{12} x_{23}} - \frac{x_{12}}{x_{13} x_{23}} - \frac{x_{23}}{x_{12} x_{13}}\biggr) }
\end{array} } \right\} 
\end{equation}
where $n = 1+k_2+k_3$. 

\noindent (2) The scalar curvature $S_g$ is given by
{ \begin{eqnarray}\label{scalarSt2}
S_g &=& \frac{1}{2}\left(\frac{d_2}{x_2} +\frac{d_{12}}{x_{12}} + \frac{d_{13}}{x_{13}} + \frac{d_{23}}{x_{23}}\right) -\frac{1}{4x_2}\left(A_{2(12)(12)} + A_{2(13)(13)} + A_{2(23)(23)}\right) \nonumber\\
&&-\frac{A_{(12)(13)(23)}}{2}\left(\frac{x_{12}}{x_{13}x_{23}} + \frac{x_{13}}{x_{12}x_{23}} + \frac{x_{23}}{x_{12}x_{13}}\right)
\end{eqnarray} }
where $d_2 = \dim\fr{m}_2$ and $d_{ij} = \dim\fr{m}_{ij}$, $i\neq j = 1,2,3$.
\end{lemma}

We normalize the metric $g = (x_2, x_{12}, x_{13}, x_{23})$ by setting $x_{23} = 1$, then by solving the system $\{r_2-r_{12}=0, r_{12}-r_{23}=0, r_{13}-r_{23} = 0\}$ for $k_2 = 4$ and $k_3 = 2$, we take the following:
\begin{theorem}\textnormal{(\cite{ArSaSt1})}
The Stiefel manifold $V_5\mathbb{R}^7 = \SO(7)/\SO(2)$  admits at least four invariant Einstein metrics, which are determined by the $\Ad(\SO(4)\times\SO(2))$-invariant inner products of the form (\ref{metric2}) given as:
$g_1 = (1.27429, 1.27429, 1, 1)$,\ $g_2 = (0.392375, 0.392375, 1, 1)$,
$g_3 = (0.245146, 1.01652, 0.253386, 1)$, and  $g_4 = (0.291175, 0.669071, 1.16137, 1)$.
\end{theorem}

\section{Generalized flag manifolds}

\subsection{Description of flag manifolds in terms of painted Dynkin diagrams} 
Let $\fr{g}$ and $\fr{k}$ be the Lie algebras of $G$ and $K$ respectively and $\fr{g}^{\bb{C}}$, $\fr{k}^{\bb{C}}$ be their complexifications.  We choose a maximal torus  $T$ in $G$ and let $\fr{h}$ be the Lie algebra of $T$.  Then the complexification $\fr{h}^{\bb{C}}$ is a Cartan subalgebra of $\fr{g}^{\bb{C}}$.  Let $R\subset(\fr{h}^{\bb{C}})^{*}$   be the root system of $\fr{g}^{\mathbb{C}}$ relative to the Cartan subalgebra $\fr{h}^{\mathbb{C}}$  and consider the  root space decomposition $\fr{g}^{\mathbb{C}}=\fr{h}^{\mathbb{C}}\oplus\sum_{\al\in R}\fr{g}_{\al}^{\mathbb{C}}$, where $\fr{g}_{\al}^{\bb{C}}=\{X\in\fr{g}^{\bb{C}} : \ad(H)X=\al(H)X,   \ \mbox{for all} \   H\in\fr{h}^{\bb{C}}\}$  denotes the root space associated to a root $\al$.  Assume that $\fr{g}^{\bb{C}}$ is semisimple, so the Killing form $B$ of $\fr{g}^{\bb{C}}$ is non degenerate, and we establish a natural isomorphism between $\fr{h}^{\bb{C}}$ and the dual space $(\fr{h}^{\bb{C}})^{*}$ as follows: for every $\al\in(\fr{h}^{\bb{C}})^{*}$ we define $H_{\al}\in\fr{h}^{\bb{C}}$ by the equation $B(H, H_{\al}) = \al(H)$, for all $H\in\fr{h}^{\bb{C}}$.  Choose a  basis $\Pi=\{\al_{1}, \ldots, \al_{\ell}\}$ $(\dim\fr{h}^{\mathbb{C}}=\ell)$ of simple roots for $R$, and let $R^{+}$ be a choise of positive roots.  

Since $\fr{h}^{\bb{C}}\subset\fr{k}^{\bb{C}}\subset\fr{g}^{\bb{C}}$,  there is a closed subsystem $R_{K}$ of $R$ such that 
$\fr{k}^{\bb{C}}=\fr{h}^{\bb{C}}\oplus\sum_{\al\in R_{K}}\fr{g}_{\al}^{\mathbb{C}}$.  In particular, we can always  find  a subset $\Pi_{K}\subset\Pi$   such that  
$
R_{K}=R\cap\left\langle\Pi_{K} \right\rangle=\{\be\in R : \be=\sum_{\al_{i}\in\Pi_{K}}k_{i}\al_{i}, \ k_{i}\in\bb{Z}\},
$
where    $\left\langle\Pi_{K} \right\rangle$  is the space of roots generated by $\Pi_{K}$  with integer coefficients. The complex Lie algebra $\fr{k}^{\mathbb{C}}$  is a maximal reductive subalgebra of $\fr{g}^{\mathbb{C}}$  
and thus it  admits the decomposition $\fr{k}^{\mathbb{C}}=\fr{z}(\fr{k}^{\mathbb{C}})\oplus\fr{k}^{\mathbb{C}}_{ss}$,   
where $\fr{z}(\fr{k}^{\mathbb{C}})$ is  the  center of $\fr{k}^{\mathbb{C}}$ and  $\fr{k}^{\mathbb{C}}_{ss}=[\fr{k}^{\mathbb{C}}, \fr{k}^{\mathbb{C}}]$ 
is its  semisimple part.  Note that $\fr{k}^{\bb{C}}_{ss}$ is given by
$\fr{k}^{\bb{C}}_{ss}=\fr{h}'\oplus\sum_{\al\in R_{K}}\fr{g}_{\al}^{\bb{C}},$
where $\fr{h}'=\sum_{\al\in\Pi_{K}}\bb{C}H_{\al}\subset\fr{h}^{\bb{C}}$ is a Cartan subalgebra of $\fr{k}^{\mathbb{C}}_{ss}$.  In fact, $R_{K}$ is  the root system of the semisimple part  $\fr{k}^{\mathbb{C}}_{ss}$ and  $\Pi_{K}$ is a corresponding basis.  Thus we easily conclude that $\dim_{\bb{C}}\fr{h}'={\rm card}\,{\Pi_{K}}$, where  ${\rm card}\,{\Pi_{K}}$ denotes  the cardinality of the set $\Pi_{K}$.
Let $K$ be the connected Lie subgroup of $G$ generated by $\fr{k}=\fr{k}^{\bb{C}}\cap \fr{g}$.  Then the homogeneous manifold $M=G/K$ is a flag manifold, and any flag manifold is defined in this way, i.e. by the choise of a triple $(\fr{g}^{\bb{c}}, \Pi, \Pi_{K})$.

Set  $\Pi_{M}=\Pi\backslash \Pi_{K}$, and  $R_{M}=R\backslash R_{K}$, such that $\Pi=\Pi_{K}\cup \Pi_{M}$, and $R=R_{K}\cup R_{M}$, respectively.  Roots in $R_{M}$ are called  {\it complementary roots}, and they play an important role in the geometry of $M=G/K$.  
For example, let $\fr{m}$ the orthogonal complement of $\fr{k}$ in $\fr{g}$ with respect to $B$.  Then we have $[\fr{k}, \fr{m}]\subset\fr{m}$ where $\fr{m}\cong T_{o}(G/K)$.  We set $R_{M}^{+} = R^{+}\backslash R_{K}^{+}$ where $R_{K}^{+}$ is the system of positive roots of $\fr{k}^{\bb{C}}$ $(R_{K}^{+}\subset R^{+})$.  Then $\fr{m}^{\mathbb{C}}=\sum_{\al\in R_{M}}\fr{g}_{\al}^{\bb{C}}$. 

We conclude that all information contained in  $\Pi=\Pi_{K}\cup\Pi_{M}$ can be presented graphically by the  painted Dynkin diagram of $M=G/K$.
\begin{definition}\label{pdd}
Let $\Gamma=\Gamma(\Pi)$ be the Dynkin diagram of the fundamental system $\Pi$.  By painting in  black the nodes 
of $\Gamma$  corresponding  to  $\Pi_{M}$, we obtain the painted Dynkin diagram of the flag manifold $G/K$. In this diagram
the subsystem $\Pi_{K}$ is determined as the subdiagram of white roots. 
\end{definition}
  
Conversely, given a painted Dynkin diagram, in order to obtain the corresponding flag manifold $M=G/K$ we are working as follows: We define $G$ as the unique simply connected Lie group  corresponding to the  underlying Dynkin diagram $\Gamma=\Gamma(\Pi)$. The connected Lie subgroup $K\subset G$ is defined by using  the additional information $\Pi=\Pi_{K}\cup\Pi_{M}$ encoded into the painted Dynkin diagram.  The semisimple part of $K$ is obtained from  the (not  necessarily connected) subdiagram of white roots, and   each black root, i.e.   each  root in $\Pi_{M}$,  gives rise to one $U(1)$-summand.  Thus the painted Dynkin diagram determines  the  isotropy subgroup $K$ and the space $M=G/K$ completely.   By using certain rules to determine whether  different painted Dynkin diagrams define isomorphic flag manifolds (see \cite{AlAr}),  one  can obtain all flag manifolds $G/K$ of  a compact  simple Lie group $G$.

\subsection{$\fr{t}$-roots and isotropy summands}
We study the isotropy representation of a generalized flag manifold $M=G/K$ of a compact simple Lie group $G$  in terms of $\fr{t}$-roots.  In order to realise the decomposition of $\fr{m}$ into irreducible $\Ad(K)$-modules we use the center $\fr{t}$ of the real Lie algebra $\fr{k}$.  For simplicity, we fix  a system of simple roots $\Pi=\{\al_1,\ldots, \al_r, \phi_1, \ldots, \phi_k\}$ of $R$,  such that $r+k=\ell=\rnk\fr{g}^{\bb{C}}$  and we assume that $\Pi_{K}=\{\phi_1, \ldots, \phi_k\}$ is a basis of the root system $R_{K}$ of $K$ so $\Pi_{M}=\Pi\backslash \Pi_{K}=\{\al_{1}, \ldots, \al_{r}\}$.  Let $\Lambda_{1}, \ldots, \Lambda_{r}$ be  the fundamental weights  corresponding to the simple roots of $\Pi_{M}$, i.e. the linear forms defined by 
$\frac{2(\Lambda_{i}, \al_{j})}{(\al_{j}, \al_{j})}=\delta_{ij}, (\Lambda_{j}, \phi_{i})=0,$
where $(\al, \be)$ denotes the inner product on $(\fr{h}^{\mathbb{C}})^{*}$ given by $(\al, \be)=(H_{\al}, H_{\be})$, for all $\al, \be\in(\fr{h}^{\mathbb{C}})^{*}$.  Then the $\{\Lambda_{i} : 1\leq i\leq r\}$ is a basis of the dual space $\fr{t}^{*}$ of $\fr{t}$,  $\fr{t}^{*}=\sum_{i=1}^{r}\mathbb{R}\Lambda_{i}$ and $\dim\fr{t}^*=\dim \fr{t}=r$.

Consider now the linear restriction map  $\kappa : \fr{h}^{*}\to \fr{t}^{*}$ defined by $\kappa(\al)=\al|_{\fr{t}}$,
and set  $R_{\fr{t}} = \kappa(R)=\kappa(R_{M})$.  
\begin{definition}
The elements of $R_{\fr{t}}$ are called {\it $\fr{t}$-roots}.  
\end{definition}

As we saw the flag manifolds $G/K$ are determined by pairs $(\fr{g}, \Pi, \Pi_K)$.  The number of $\ad(\fr{k})$-submodules of $\fr{m}\cong T_{o}(G/K)$ correspond to the {\it Dynkin mark} of the simple root we paint black on the Dynkin diagram.  We recall the following definition

\begin{definition}\label{mark}
The Dynkin mark of a simple root $\al_{i}\in\Pi$  $(i=1, \ldots, \ell)$, is  the  positive integer $m_{i}$ in 
the expression of the highest root $\widetilde{\al}=\sum_{k=1}^{\ell}m_{k}\al_{k}$ in terms of simple roots.  We will denote by $\Mark$ the function $\Mark : \Pi\to\bb{Z}^{+}$  with $\Mark(\al_{i})=m_{i}$.
\end{definition}

A fundamental result about $\fr{t}$-root is the following:

\begin{prop}\textnormal{(\cite{AlPe})} \label{isotropy} 
There exists a one-to-one correspondence between $\fr{t}$-roots $\xi$ and irreducible $\ad(\fr{k}^{\mathbb{C}})$-submodules $\fr{m}_{\xi}$\footnote{We mean that $[\fr{k}^{\bb{C}}, \fr{m}_{\xi}]\subset\fr{m}_{\xi}$ for all $\xi\in R_{\fr{t}}$.} of the isotropy representation of $\fr{m}^{\mathbb{C}}$, which is given by   
$
R_{\fr{t}}\ni\xi \ \leftrightarrow \ \fr{m}_{\xi} =\sum_{\al\in R_{M}: \kappa(\al)=\xi}\mathbb{C}E_{\al}.
$
Thus $\fr{m}^{\mathbb{C}} = \bigoplus_{\xi\in R_{\fr{t}}} \fr{m}_{\xi}.$  Moreover, these submodules are non equivalent as $\ad(\fr{k}^{\mathbb{C}})$-modules.		 		   
\end{prop}


\subsection{Flag manifolds with four isotropy summands}

The generalized flag manifolds with four isotropy summands can be separated into two types I and II.  Type I is defined by the set: $\Pi\backslash \Pi_K = \{\al_i : \Mark(\al_{i})=4\}$ (that is $b_2(G/K) = 1$) and Type II is given by $\Pi\backslash \Pi_K = \{\al_i, \al_j : \Mark(\al_{i})=1, \Mark(\al_{i})=2\}$ (however this set may be define flag manifolds with four or five isotropy summands).  The classification of those spaces was given by A. Arvanitoyeorgos and I. Chrysikos in \cite{ArCh2}.  Below we give the flag manifolds together with the Dynkin diagram of Type I.

\begin{eqnarray*}
{\F_4/\SU(3)\times \SU(2)\times \U(1)} \ \ \ \ \ \ \  && \ \ \ \quad \ \ \ \quad \E_7/\SU(4)\times \SU(3)\times \SU(2)\times \U(1) \\
\begin{picture}(160,20)(58, 5)
\put(87,10){\circle{4 }}
\put(87,18){\makebox(0,0){$\al_1$}}
\put(87,2){\makebox(0,0){2}}
\put(89,10.5){\line(1,0){14}}
\put(105,10){\circle{4 }}
\put(105,18){\makebox(0,0){$\al_2$}}
\put(105,2){\makebox(0,0){3}}
\put(107, 11.3){\line(1,0){12.4}}
\put(107, 9.1){\line(1,0){12.6}}
\put(117.6, 8){\scriptsize $>$}
\put(124.5,10){\circle*{4 }}
\put(124.5, 18){\makebox(0,0){$\al_3$}}
\put(124.5, 2){\makebox(0,0){4}}
\put(126,10.5){\line(1,0){16}}
\put(144,10){\circle{4 }}
\put(144,18){\makebox(0,0){$\al_4$}}
\put(144,2){\makebox(0,0){2}}
\end{picture}
&& 
\begin{picture}(160,25)(-32, 6)
\put(15, 9.5){\circle{4 }}
\put(15, 18){\makebox(0,0){$\al_1$}}
\put(15,2){\makebox(0,0){1}}
\put(17, 10){\line(1,0){14}}
\put(33.5, 9.5){\circle{4 }}
\put(33.5, 18){\makebox(0,0){$\al_2$}}
\put(33.5,2){\makebox(0,0){2}}
\put(35, 10){\line(1,0){13.6}}
 \put(51, 9.5){\circle{4 }}
 \put(51, 18){\makebox(0,0){$\al_3$}}
\put(51,2){\makebox(0,0){3}}
\put(69,-6){\line(0,1){14}}
\put(53,10){\line(1,0){14}}
\put(69,9.5){\circle*{4 }}
\put(69.9, 18){\makebox(0,0){$\al_4$}}
\put(65,2){\makebox(0,0){4}}
\put(69,-8){\circle{4}}
\put(79, -10){\makebox(0,0){$\al_7$}}
\put(62,-10){\makebox(0,0){2}}
\put(71,10){\line(1,0){16}}
\put(89,9.5){\circle{4 }}
\put(89, 18){\makebox(0,0){$\al_5$}}
\put(89,2){\makebox(0,0){3}}
\put(90.7,10){\line(1,0){16}}
\put(109,9.5){\circle{4 }}
\put(109, 18){\makebox(0,0){$\al_6$}}
\put(109, 2){\makebox(0,0){2}}
\end{picture}
\\ \\
{\E_{8}(\al_3)}/\SO(10)\times \SU(3)\times \U(1) && \ \ \ \quad \ \ \ \quad  {\E_{8}(\al_6)}/\SU(7)\times \SU(2)\times \U(1) \\
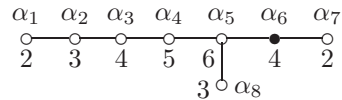
\begin{picture}(160,25)(-14, 6)
\put(15, 9.5){\circle{4 }}
\put(15, 18){\makebox(0,0){$\al_1$}}
\put(15,2){\makebox(0,0){2}}
\put(17, 10){\line(1,0){14}}
\put(33.5, 9.5){\circle{4 }}
\put(33.5, 18){\makebox(0,0){$\al_2$}}
\put(33.5,2){\makebox(0,0){3}}
\put(35, 10){\line(1,0){13.6}}
 \put(51, 9.5){\circle*{4 }}
 \put(51, 18){\makebox(0,0){$\al_3$}}
\put(51,2){\makebox(0,0){4}}
\put(53,10){\line(1,0){14}}
\put(69,9.5){\circle{4 }}
\put(69, 18){\makebox(0,0){$\al_4$}}
\put(69,2){\makebox(0,0){5}}
\put(89,-8){\circle{4}}
\put(99, -9.5){\makebox(0,0){$\al_8$}}
\put(82,-9.5){\makebox(0,0){3}}
\put(89,-6){\line(0,1){14}}
\put(71,10){\line(1,0){16}}
\put(89,9.5){\circle{4 }}
\put(89, 18){\makebox(0,0){$\al_5$}}
\put(84,2){\makebox(0,0){6}}
\put(90.7,10){\line(1,0){16}}
\put(109,9.5){\circle{4 }}
\put(109, 18){\makebox(0,0){$\al_6$}}
\put(109,2){\makebox(0,0){4}}
\put(111,10){\line(1,0){16}}
\put(129,9.5){\circle{4 }}
\put(129, 18){\makebox(0,0){$\al_7$}}
\put(129,2){\makebox(0,0){2}}
\end{picture}
&&
\begin{picture}(160,25)(-32, 6)
\put(15, 9.5){\circle{4 }}
\put(15, 18){\makebox(0,0){$\al_1$}}
\put(15,2){\makebox(0,0){2}}
\put(17, 10){\line(1,0){14}}
\put(33.5, 9.5){\circle{4 }}
\put(33.5, 18){\makebox(0,0){$\al_2$}}
\put(33.5,2){\makebox(0,0){3}}
\put(35, 10){\line(1,0){13.6}}
 \put(51, 9.5){\circle{4 }}
 \put(51, 18){\makebox(0,0){$\al_3$}}
\put(51,2){\makebox(0,0){4}}
\put(53,10){\line(1,0){14}}
\put(69,9.5){\circle{4 }}
\put(69, 18){\makebox(0,0){$\al_4$}}
\put(69,2){\makebox(0,0){5}}
\put(89,-8){\circle{4}}
\put(99, -9.5){\makebox(0,0){$\al_8$}}
\put(82,-9.5){\makebox(0,0){3}}
\put(89,-6){\line(0,1){14}}
\put(71,10){\line(1,0){16}}
\put(89,9.5){\circle{4 }}
\put(89, 18){\makebox(0,0){$\al_5$}}
\put(84,2){\makebox(0,0){6}}
\put(90.7,10){\line(1,0){16}}
\put(109,9.5){\circle*{4 }}
\put(109, 18){\makebox(0,0){$\al_6$}}
\put(109,2){\makebox(0,0){4}}
\put(111,10){\line(1,0){16}}
\put(129,9.5){\circle{4 }}
\put(129, 18){\makebox(0,0){$\al_7$}}
\put(129,2){\makebox(0,0){2}}
\end{picture} 
\\
\end{eqnarray*}

Let $M = G/K$ be a generalized flag manifold of Type I and let $\fr{m} = \fr{m}_1\oplus\fr{m}_2\oplus\fr{m}_3\oplus\fr{m}_4$ be a decomposition of $\fr{m} = T_{o}M$ into irreducible non-equivalent $\Ad(K)$-modules, with respect to the negative of the Killing form $B$ of $G$. Then, a $G$-invariant metric on $M = G/K$ is given by
\begin{equation}\label{metricFlag}
g = \langle\cdot, \cdot\rangle = x_1(−B)|_{\fr{m}_1} + x_2(−B)|_{\fr{m}_2} + x_3(−B)|_{\fr{m}_3} + x_4(−B)|_{\fr{m}_4} 
\end{equation}
where $x_i\in\bb{R}_{+}$, $i=1,2,3,4$. Very often we will denote such metrics with $g = (x_1, x_2, x_3, x_4)$.  The Ricci tensor for the above metric has (as symmetric, covariant 2-tensor) the same expression that is:
$
\Ric_{\langle\cdot, \cdot\rangle} =\sum_{i=1}^4 x_i r_i(-B)|_{\fr{m}_i}, 
$
where $r_i, i=1,2,3,4$ are the components of the Ricci tensor and are given as follows:

\begin{prop}\textnormal{(\cite{ArCh2})}
(1) The components  of  the Ricci tensor $\Ric$ for the invariant metric $\langle\cdot, \cdot\rangle$ on $G/K$ defined by {\em{(\ref{metricFlag})}}, are given as follows:  
 \begin{equation}\label{compI}
\left. {\small \begin{array}{l} 
r_1=\displaystyle\frac{1}{2x_1}-\frac{A_{112}}{2d_1}\frac{x_2}{x_1^2}+
\frac{A_{123}}{2d_1}\Big(\frac{x_1}{x_2x_3}-\frac{x_2}{x_1x_3}-\frac{x_3}{x_1x_2}\Big)+
\frac{A_{134}}{2d_1}\Big(\frac{x_1}{x_3x_4}-\frac{x_3}{x_1x_4}-\frac{x_4}{x_1x_3}\Big) \\ \\
r_2=\displaystyle\frac{1}{2x_2}-\frac{A_{224}}{2d_2}\frac{x_4}{x_2^2}+
\frac{A_{112}}{4d_2}\Big(\frac{x_2}{x_1^2}-\frac{2}{x_2}\Big)+
\frac{A_{123}}{2d_2}\Big(\frac{x_2}{x_1x_3}-\frac{x_1}{x_2x_3}-\frac{x_3}{x_1x_2}\Big) \\ \\
r_3=\displaystyle\frac{1}{2x_3}+
\frac{A_{123}}{2d_3}\Big(\frac{x_3}{x_1x_2}-\frac{x_2}{x_1x_3}-\frac{x_1}{x_2x_3}\Big)+
\frac{A_{134}}{2d_3}\Big(\frac{x_3}{x_1x_4}-\frac{x_1}{x_3x_4}-\frac{x_4}{x_1x_3}\Big) \\ \\ 
r_4=\displaystyle\frac{1}{2x_4}+
\frac{A_{224}}{4d_4}\Big(\frac{x_4}{x_2^2}-\frac{2}{x_4}\Big)+
\frac{A_{134}}{2d_4}\Big(\frac{x_4}{x_1x_3}-\frac{x_1}{x_3x_4}-\frac{x_3}{x_1x_4}\Big). 
\end{array} } \right\} 
\end{equation}

\noindent(2) The scalar curvature is given by
\begin{eqnarray}\label{scalFlag}
S_g &=& \frac{1}{2}\left(\frac{d_1}{x_1} + \frac{d_2}{x_2} + \frac{d_3}{x_3} + \frac{d_4}{x_4}\right) + \frac{A_{112}}{4}\left(-\frac{2}{x_2}-\frac{x_2}{x_1^2}\right) + \frac{A_{123}}{2}\left(-\frac{x_{1}}{x_2x_3}-\frac{x_2}{x_1x_3}-\frac{x_3}{x_1x_2}\right) \nonumber\\
&& + \frac{A_{224}}{4}\left(-\frac{2}{x_4}-\frac{x_4}{x_2^2}\right)+ \frac{A_{134}}{2}\left(-\frac{x_1}{x_3x_4}-\frac{x_3}{x_1x_4}-\frac{x_4}{x_1x_3}\right),  
\end{eqnarray}
where $d_i = \dim\fr{m}_i$, $i=1,2,3,4$.
\end{prop}
In \cite{ArCh2} the authors computed the numbers $A_{ijk}$ using the twistor fibration which admits any flag manifold $M = G/K$ of a compact (semi)-simple Lie group $G$, over an irreducible symmetric space $G/L$ of compact type.  
In particular we have the following table for $A_{ijk}$ and the dimensions $d_i$.

\smallskip
 \begin{center}
{\small {\bf Table 4.} The numbers $A_{ijk}$ and the dimensions for the flag manifold of Type I}  
\end{center} 
{\small\small  \begin{center}
\begin{tabular}{  lllll||llll  }
\thickline
 $M=G/K$ & $A_{224}$ & $A_{112}$ & $A_{123}$ & $A_{134}$ & $d_1$ & $d_2$ & $d_3$ & $d_4$\\
\thickline
$\F_4/\SU(3)\times \SU(2)\times \U(1)$ & $2$ & $2$ & $1$ & $2/3$ & $12$ & $18$ & $4$ & $6$\\
$\E_7/\SU(4)\times \SU(3)\times \SU(2)\times \U(1)$ & $2$ & $8$ & $4$ & $4/3$ & $48$ & $36$ & $16$ & $6$\\
$\E_8(\al_3)/\SO(10)\times \SU(3)\times \U(1)$ & $2$ & $16$ & $8$ & $8/5$ & $96$ & $60$ & $32$ & $6$\\
$\E_8(\al_6)/\SU(7)\times \SU(2)\times \U(1)$ & $14/3$ & $14$ & $7$ & $14/5$ & $84$ & $70$ & $28$ & $14$\\
\thickline
\end{tabular}
\end{center}}
\medskip
\smallskip

We normalize the metric $g = (x_1, x_2, x_3, x_4)$ by setting $x_1 = 1$.  Then $g$ is Einstein if and only the system: $\{r_1-r_2=0,\ r_2-r_3=0, \ r_3-r_4=0\}$ has positive solution.  After solving the previous system for each flag manifold separately we obtain the following theorem

\begin{theorem}\textnormal{(\cite{ArCh2})}\label{TheoremFlag}
(1) The generalized flag manifolds $G/K$ associated to the exceptional Lie groups $\F_4, \E_7$ and $\E_8(\al_3)$ admits (up to scale) three $G$-invariant Einstein metrics.  One is K\"ahler given by $g = (1,2,3,4)$ and other two are non-K\"ahler given approximatelly as follows:
\begin{eqnarray*}
\F_4 &:& g_1= (1, 1.2761, 1.9578, 2.3178),\ g_2= (1, 0.9704, 0.2291, 1.0097)\\
\E_7 &:& g_1= (1, 0.8233, 1.2942, 1.3449),\ g_2= (1, 0.9912, 0.5783, 1.1312)\\
\E_8(\al_3) &:& g_1= (1, 0.9133, 1.4136, 1.5196),\ g_2= (1, 0.9663, 0.4898, 1.0809)
\end{eqnarray*}
(2) If $G/K = \E_8(\al_6)/\SO(10)\times\SU(3)\times\U(1)$ then $G/K$ admits (up to scale) five $E_8$-invariant Einstein metrics. One is K\"ahler given by $g = (1,2,3,4)$ and other four are non-K\"ahler given approximatelly as follows: 
$
g_1 = (1, 0.6496, 1.1094, 1.0610),  g_2 = (1, 1.1560, 1.0178, 0.2146), 
g_3 = (1, 1.0970, 0.7703, 1.2969),  g_4 = (1, 0.7633, 1.0090, 0.1910).
$
\end{theorem}

\section{Poincar\'e compactification}

The method of Poincar\'e compactification dates back to 1881.  The main idea is to pass the study of a vector field on a non compact manifold, to its study on the sphere (compact manifold).   This allows us to better understand its behavior at infinity.  
Poincar\'e was studying the behavior of polynomial planar vector fields at infinity by means of the central projection.  For more details of the description of this method for $n$-dimensional case can be found in \cite{Ve}.  Next, we will describe the method in three dimensions.

Let $(x_1, x_2, x_3)$ be the coordinates of $\bb{R}^3$ and $X = P_{1}(x_1, x_2, x_3)\frac{\partial}{\partial x_1}  +  P_{2}(x_1, x_2, x_3)\frac{\partial}{\partial x_2} + P_3(x_1, x_2, x_3)\frac{\partial}{\partial x_3}$ be a polynomial vector field of degree $d = {\rm max}\{{\rm deg}(P_1), {\rm deg}(P_2), {\rm deg}(P_3)\}$.  We consider the sphere $\bb{S}^3 = \{(y_1, y_2, y_3, y_4)\in\bb{R}^4 : y_1^2 + y_2^2 + y_3^2 + y_4^2 = 1\}$, which we shall call the Poincar\'e sphere with north hemisphere $\bb{S}^3_{+} = \{y\in \bb{S}^3 : y_4>0\}$, south hemishpere $\bb{S}^3_{-} = \{y\in\bb{S}^3 : y_4<0\}$ and equator $\bb{S}^2 = \{y\in\bb{S}^3 : y_4 = 0\}$.

The central projection from $\bb{R}^3$ to the Poincar\'e sphere is defined as follows: 
\begin{eqnarray}\label{projection}
f_{+} : \bb{R}^3 \to \bb{S}^3, && (x_1, x_2, x_3) \mapsto \left(\frac{x_1}{\Delta(x)}, \frac{x_2}{\Delta(x)}, \frac{x_3}{\Delta(x)}, \frac{1}{\Delta(x)}\right) \nonumber\\
f_{-} : \bb{R}^3 \to \bb{S}^3, && (x_1, x_2, x_3) \mapsto \left(\frac{-x_1}{\Delta(x)}, \frac{-x_2}{\Delta(x)}, \frac{-x_3}{\Delta(x)}, \frac{-1}{\Delta(x)}\right), 
\end{eqnarray}  
where $\Delta(x) = \sqrt{x_1^2 + x_2^2 + x_3^2 + 1}$.  The maps $f_{+}$ and $f_{-}$ define the following two vectors fields on each hemisphere
$$
(df_{+})_xX(x) \ \mbox{at}\ y = f_{+}(x), \ \mbox{and} \ (df_{-})_xX(x) \ \mbox{at}\ y = f_{-}(x).
$$
However, these two vector fields can not be extended to the equator $\bb{S}^2$.  Indeed, in the case where $y\in U_1 = \{y\in\bb{S}^3 : y_1 > 0\}$ where $U_1$ is the local chart with corresponding map: 
\begin{equation}\label{chart1}
\varphi_1 : U_1 \to \bb{R}^3, \ \varphi_1(y) = \left(\frac{y_2}{y_1}, \frac{y_3}{y_1}, \frac{y_4}{y_1}\right)
\end{equation}
we shall denote by $z = (z_1, z_2, z_3)$ the value of $\varphi_1(y)$, so  that $z$ represent different things according to the case under consideration.

\begin{remark}
\textnormal{There are eight of these charts: $U_i = \{y\in\bb{S}^3 : y_i > 0\}$ and $V_i = \{y\in\bb{S}^3 : y_i < 0\}$ with $i = 1,2,3,4$ and the local maps for the corresponding charts are given by
$
\varphi_i : U_i \to \bb{R}^3, \psi_i : V_i \to \bb{R}^3
$
with
$$
\varphi_i(y) = -\psi_i(y) = \left(\frac{y_n}{y_i}, \frac{y_m}{y_i}, \frac{y_k}{y_i}\right), \ \mbox{for} \ n<m<k \ \mbox{and}\ n,m,k\neq i.
$$
}
\end{remark}

By using (\ref{projection}) and (\ref{chart1}) we have
\begin{equation}\label{astro}
z = \varphi_1(y) = \left(\frac{y_2}{y_1}, \frac{y_3}{y_1}, \frac{y_4}{y_1}\right) = \left(\frac{x_2}{x_1}, \frac{x_3}{x_1}, \frac{1}{x_1}\right).
\end{equation}
The differential of $f_{+}$ in the case where $y$ is in the upper hemisphere, is given by the following matrix 
$$
\{(df_{+})_x\}_{ij} ={\frac{\partial z_i}{\partial x_j}} = \begin{pmatrix}
-\frac{x_2}{x_1^2} & \frac{1}{x_1} & 0\\ 
-\frac{x_3}{x_1^2} & 0 & \frac{1}{x_1}\\
-\frac{1}{x_1^2} & 0 & 0
\end{pmatrix}.
$$
Therefore
\begin{equation}\label{diaforiko1}
(df_{+})_xX(x) = (1/x_1)^2 \left(-x_2 P_1(x) + x_2P_2(x), -x_3P_1(x)+x_1P_3(x), -P_1(x)\right).
\end{equation}
We can write the above expression in terms of the corresponding point on the sphere and we get for (\ref{astro})
\begin{eqnarray}\label{diaforiko2}
(df_{+})_xX(x) & = & \left(\frac{y_4}{y_1}\right)^2 \left\lbrace-\frac{y_2}{y_4}P_1\left(\frac{y_1}{y_4}, \frac{y_2}{y_4}, \frac{y_3}{y_4}\right) + \frac{y_1}{y_4}P_2\left(\frac{y_1}{y_4}, \frac{y_2}{y_4}, \frac{y_3}{y_4}\right),  \right. 
\nonumber\\
& & \quad \quad \quad \left.  -\frac{y_3}{y_4}P_1\left(\frac{y_1}{y_4}, \frac{y_2}{y_4}, \frac{y_3}{y_4}\right) + \frac{y_1}{y_4}P_3\left(\frac{y_1}{y_4}, \frac{y_2}{y_4}, \frac{y_3}{y_4}\right), -P_1\left(\frac{y_1}{y_4}, \, \frac{y_2}{y_4}, \frac{y_3}{y_4}\right)        \right\rbrace 
\end{eqnarray}

It is clear that if we multiply the above vector field by the factor $y_{4}^{d-1}$ (it depends only on the degree $d$ of the polynomial vector field $X$) then the vector field (\ref{diaforiko2}) extends into $\bb{S}^2$.
We will take the same expression for the induced vector field as (\ref{diaforiko2}) in case where we work on lower hemisphere.  Actually, if we denote by $\bar{X}(y)$ the vector field on $\bb{S}^3\setminus\bb{S}^2 = \bb{S}^3_{+}\cup\bb{S}^3_{-}$ then to extend $\bar{X}(y)$ to the Poincar\'e sphere $\bb{S}^3$, we define the Poincar\'e compactification of $X = (P_1, P_2, P_3)$ denoted by $p(X)$ and is given as $p(X)(y) = y_4^{d-1}\bar{X}(y)$.    

\begin{theorem}\textnormal{(\cite{Ve})}
The vector field $p(X)$ extends $\bar{X}(y)$ analytically to the whole sphere $\bb{S}^3$, in such a way that the equator $\bb{S}^2$ is invariant.
\end{theorem}

The important point here is the fact that if we know the behavior of $p(X)$ around the equator, then we know the behavior of $X$ in the neighborhood of infinity.

Now we will give explicitly the expressions of $p(X)$ in the local charts.  It is convenient to express these fields in terms of the variable $z = (z_1, z_2, z_3)$ of $\bb{R}^3$.  In the chart $(U_1, \varphi_1)$ we have from (\ref{astro}) that $y_1/y_4 = 1/z_3$, $y_3/y_4 = z_2/z_3$, $y_2/y_4 = z_1/z_3$ and also $y_4^2 = z_3^2/\Delta(z)^2$ where $\Delta(z) = \sqrt{z_1^2 + z_2^2 + z_3^2 + 1}$.  From (\ref{astro}) and because $y_1>0$ the $y_4$ has the same sign as $z_3$.  Thus $y_4 = z_3/\Delta(z)$.  Then after some substitutions it turns out that $p(X)$ can be expressed as 
\begin{eqnarray}\label{Xartis1}
&& \frac{z^{d}_3}{\Delta(z)^{d-1}}\left\lbrace -z_1 P_{1}\left(\frac{1}{z_3}, \frac{z_1}{z_3}, \frac{z_2}{z_3}\right) + P_2\left(\frac{1}{z_3}, \frac{z_1}{z_3}, \frac{z_2}{z_3}\right),  \right. \nonumber\\
&& \quad \quad \quad \quad \ \ \left. -z_2 P_{1}\left(\frac{1}{z_3}, \frac{z_1}{z_3}, \frac{z_2}{z_3}\right) + P_3\left(\frac{1}{z_3}, \frac{z_1}{z_3}, \frac{z_2}{z_3}\right), -z_3  P_{1}\left(\frac{1}{z_3}, \frac{z_1}{z_3}, \frac{z_2}{z_3}\right)  \right\rbrace.
\end{eqnarray} 

The expression in the chart $(U_2, \varphi_2)$ is
\begin{eqnarray}\label{Xartis2}
&& \frac{z_3^{d}}{\Delta(z)^{d-1}}\left\lbrace P_{1}\left(\frac{z_1}{z_3}, \frac{1}{z_3}, \frac{z_2}{z_3}\right) - z_1 P_2\left(\frac{z_1}{z_3}, \frac{1}{z_3}, \frac{z_2}{z_3}\right),  \right. \nonumber\\
&& \quad \quad \quad \quad \left. -z_2 P_{2}\left(\frac{z_1}{z_3}, \frac{1}{z_3}, \frac{z_2}{z_3}\right) + P_3\left(\frac{z_1}{z_3}, \frac{1}{z_3}, \frac{z_2}{z_3}\right), -z_3  P_{2}\left(\frac{z_1}{z_3}, \frac{1}{z_3}, \frac{z_2}{z_3}\right)  \right\rbrace.
\end{eqnarray}
    
For the chart $(U_3, \varphi_3)$ we have
\begin{eqnarray}\label{Xartis3}
&& \frac{z_3^{d}}{\Delta(z)^{d-1}}\left\lbrace P_{1}\left(\frac{z_1}{z_3}, \frac{z_2}{z_3}, \frac{1}{z_3}\right) - z_1 P_3\left(\frac{z_1}{z_3}, \frac{z_2}{z_3}, \frac{1}{z_3}\right),  \right. \nonumber\\
&& \quad \quad \quad \quad \ \ \left. P_{2}\left(\frac{z_1}{z_3}, \frac{z_2}{z_3}, \frac{1}{z_3}\right) - z_2 P_3\left(\frac{z_1}{z_3}, \frac{z_2}{z_3}, \frac{1}{z_3}\right), -z_3  P_{3}\left(\frac{z_1}{z_3}, \frac{z_2}{z_3}, \frac{1}{z_3}\right)  \right\rbrace.
\end{eqnarray}

Finally, for the chart $(U_4, \varphi_4)$ we have
$\lbrace P_{1}(z_1, z_2, z_3),  P_{2}(z_1, z_2, z_3), P_{3}(z_1, z_2, z_3)\rbrace$.
We can avoid the use of the factor $1/\Delta(z)^{d-1}$ in the expression of $p(X)$.  Also, note that for the singularities at infinity have $z_3 = 0$.  For the charts $(V_1, \psi_1), (V_2, \psi_2)$ and $(V_3, \psi_3)$ the expression of $p(X)$ is the same as (\ref{Xartis1}), (\ref{Xartis2}) and (\ref{Xartis3}) multiplied by the factor $(-1)^{d-1}$.

By the same method we can describe the necessary formulas for the compactified vector field in case of four-dimensional.  We use $(z_1, z_2, z_3, z_4)$ as coordinates.  If the original polynomial vector field is $X = \sum_{i=1}^4P_{i}(x_1, x_2, x_3, x_4)\frac{\partial}{\partial x_i}$, with $d = {\rm max}\{{\rm deg}(P_1),\ldots,{\rm deg}(P_4)\}$ the degree of $X$ then the equations of the compactified field $p(X)$ are given as follows:

On the chart $(U_1, \varphi_1)$ it is
\begin{eqnarray}\label{Xartis11}
&& \frac{z^{d}_4}{\Delta(z)^{d-1}}\left\lbrace -z_1 P_{1}\left(\frac{1}{z_4}, \frac{z_1}{z_4}, \frac{z_2}{z_4},\frac{z_3}{z_4}\right) + P_2\left(\frac{1}{z_4}, \frac{z_1}{z_4}, \frac{z_2}{z_4},\frac{z_3}{z_4}\right),  \right. \nonumber\\
&& \quad \quad \quad  \left. -z_2 P_{1}\left(\frac{1}{z_4}, \frac{z_1}{z_4}, \frac{z_2}{z_4},\frac{z_3}{z_4}\right) + P_3\left(\frac{1}{z_4}, \frac{z_1}{z_4}, \frac{z_2}{z_4},\frac{z_3}{z_4}\right), -z_3  P_{1}\left(\frac{1}{z_4}, \frac{z_1}{z_4}, \frac{z_2}{z_4},\frac{z_3}{z_4}\right) + P_4\left(\frac{1}{z_4}, \frac{z_1}{z_4}, \frac{z_2}{z_4},\frac{z_3}{z_4}\right), \right.\nonumber\\
&&\quad \quad \quad \left. -z_4 P_1\left(\frac{1}{z_4}, \frac{z_1}{z_4}, \frac{z_2}{z_4},\frac{z_3}{z_4}\right) \right\rbrace
\end{eqnarray} 

The expression in the chart $(U_2, \varphi_2)$ is
\begin{eqnarray}\label{Xartis22}
&& \frac{z_4^{d}}{\Delta(z)^{d-1}}\left\lbrace P_{1}\left(\frac{z_1}{z_4}, \frac{1}{z_4}, \frac{z_2}{z_3}, \frac{z_3}{z_4}\right) - z_1 P_2\left(\frac{z_1}{z_4}, \frac{1}{z_4}, \frac{z_2}{z_3}, \frac{z_3}{z_4}\right),  \right. \nonumber\\
&& \quad \quad \quad \left. -z_2 P_{2}\left(\frac{z_1}{z_4}, \frac{1}{z_4}, \frac{z_2}{z_3}, \frac{z_3}{z_4}\right) + P_3\left(\frac{z_1}{z_4}, \frac{1}{z_4}, \frac{z_2}{z_3}, \frac{z_3}{z_4}\right), -z_3  P_{2}\left(\frac{z_1}{z_4}, \frac{1}{z_4}, \frac{z_2}{z_3}, \frac{z_3}{z_4}\right) + P_4\left(\frac{z_1}{z_4}, \frac{1}{z_4}, \frac{z_2}{z_3}, \frac{z_3}{z_4}\right), \right.\nonumber\\
&&\quad \quad \quad \left. -z_4 P_2\left(\frac{z_1}{z_4}, \frac{1}{z_4}, \frac{z_2}{z_3}, \frac{z_3}{z_4}\right)  \right\rbrace
\end{eqnarray}
    
For the chart $(U_3, \varphi_3)$ we have
\begin{eqnarray}\label{Xartis33}
&& \frac{z_4^{d}}{\Delta(z)^{d-1}}\left\lbrace P_{1}\left(\frac{z_1}{z_4}, \frac{z_2}{z_4}, \frac{1}{z_4}, \frac{z_3}{z_4}\right) - z_1 P_3\left(\frac{z_1}{z_4}, \frac{z_2}{z_4}, \frac{1}{z_4}, \frac{z_3}{z_4}\right),  \right. \nonumber\\
&& \quad \quad \quad  \left. P_{2}\left(\frac{z_1}{z_4}, \frac{z_2}{z_4}, \frac{1}{z_4}, \frac{z_3}{z_4}\right) - z_2 P_3\left(\frac{z_1}{z_4}, \frac{z_2}{z_4}, \frac{1}{z_4}, \frac{z_3}{z_4}\right), -z_3  P_{3}\left(\frac{z_1}{z_4}, \frac{z_2}{z_4}, \frac{1}{z_4}, \frac{z_3}{z_4}\right) + P_4\left(\frac{z_1}{z_4}, \frac{z_2}{z_4}, \frac{1}{z_4}, \frac{z_3}{z_4}\right), \right.\nonumber\\
&&\quad \quad \quad   \left. -z_4 P_3\left(\frac{z_1}{z_4}, \frac{z_2}{z_4}, \frac{1}{z_4}, \frac{z_3}{z_4}\right)  \right\rbrace
\end{eqnarray}

For the chart $(U_4, \varphi_4)$ we have
\begin{eqnarray}\label{Xartis44}
&& \frac{z_4^{d}}{\Delta(z)^{d-1}}\left\lbrace P_{1}\left(\frac{z_1}{z_4}, \frac{z_2}{z_4}, \frac{z_3}{z_4}, \frac{1}{z_4}\right) - z_1 P_4\left(\frac{z_1}{z_4}, \frac{z_2}{z_4}, \frac{z_3}{z_4}, \frac{1}{z_4}\right),  \right. \nonumber\\
&& \quad \quad \quad  \left. P_{2}\left(\frac{z_1}{z_4}, \frac{z_2}{z_4}, \frac{z_3}{z_4}, \frac{1}{z_4}\right) - z_2 P_4\left(\frac{z_1}{z_4}, \frac{z_2}{z_4}, \frac{z_3}{z_4}, \frac{1}{z_4}\right), -z_3 P_{4}\left(\frac{z_1}{z_4}, \frac{z_2}{z_4}, \frac{z_3}{z_4}, \frac{1}{z_4}\right) + P_3\left(\frac{z_1}{z_4}, \frac{z_2}{z_4}, \frac{z_3}{z_4}, \frac{1}{z_4}\right), \right.\nonumber\\
&&\quad \quad \quad   \left. -z_4 P_4\left(\frac{z_1}{z_4}, \frac{z_2}{z_4}, \frac{z_3}{z_4}, \frac{1}{z_4}\right)  \right\rbrace
\end{eqnarray}

Finally for the chart $(U_5, \varphi_5)$ we have
$\lbrace P_{1}(z_1, z_2, z_3, z_4), P_{2}(z_1, z_2, z_3, z_4), P_{3}(z_1, z_2, z_3, z_4),  P_{4}(z_1, z_2, z_3, z_4) \rbrace $.  

\section{The behavior of the normalized Ricci flow}
We analyse the Ricci flow of invariant metrics on the generalized Wallach space, the Stiefel manifolds $V_2\bb{R}^n$, $V_{1+k_2}\bb{R}^{1+k_2+k_3}$ and the generalized flag manifolds $G/K$ with four isotropy summands and $b_2(G/K) = 1$.

\subsection{Ricci flow for the generalized Wallach spaces}
By using the Ricci components and scalar curvature of subsection 3.1 it is easy to see that 
system (\ref{normalRicci}) reduces to
{ \begin{eqnarray*}
\dot{x}_1 &=& \frac{1}{{x_2} {x_3}({d_1}+{d_2}+{d_3})}\big\lbrace {a_1} \left(({d_2}+{d_3})\left({x_1}^2-{x_2}^2-{x_3}^2\right)-2 {d_1} \left({x_2}^2+{x_3}^2\right)\right) \nonumber\\
&& +{x_3} (2 {d_1} {x_2}+{d_2}({x_1}+{x_2}))+{d_3} {x_2} ({x_1}+{x_3})  \big\rbrace \nonumber\\
\dot{x}_2 &=& \frac{1}{{x_1}
   {x_3} ({d_1}+{d_2}+{d_3})} \big\lbrace -{a_1} {d_1} \left({x_1}^2+{x_2}^2+{x_3}^2\right)-{a_2}
   ({d_1}+{d_2}+{d_3})\left({x_1}^2-{x_2}^2+{x_3}^2\right) \nonumber\\
&&+{x_3} ({x_1} ({d_1}+2{d_2}+{d_3})+{d_1} {x_2})+{d_3} {x_1} {x_2} \big\rbrace \nonumber\\
   \dot{x}_3 &=&\frac{1}{{x_1} {x_2} ({d_1}+{d_2}+{d_3})} \big\lbrace -{a_1} {d_1} \left({x_1}^2+{x_2}^2+{x_3}^2\right)-    {a_3}({d_1}+{d_2}+{d_3})\left({x_1}^2+{x_2}^2-{x_3}^2\right)  \nonumber\\
&& +{x_1} {x_2}({d_1}+{d_2}+2 {d_3})+{x_3} ({d_1} {x_2}+{d_2}{x_1}) \big\rbrace.
\end{eqnarray*} }
The above system is not polynomial, hence in order to apply the Poincar\'e compactification, we multiply it by the factor $(d_1+d_2+d_3)x_1x_2x_3$.  Then we obtain the following system:
\begin{eqnarray}\label{systemWa}
\dot{x}_1 &=& {x_1} \left(-2 {a_1} {d_1} {x_2}^2-2 {a_1} {d_1}{x_3}^2+{a_1} {d_2} {x_1}^2-{a_1} {d_2} {x_2}^2-{a_1}{d_2} {x_3}^2+{a_1} {d_3} {x_1}^2-{a_1} {d_3}{x_2}^2-{a_1} {d_3} {x_3}^2 \right. \nonumber\\
&&\left. +2 {d_1} {x_2} {x_3}+{d_2}{x_1} {x_3}+{d_2} {x_2} {x_3}+{d_3} {x_1}{x_2}+{d_3} {x_2} {x_3}\right) \nonumber\\
\dot{x}_2 &=& {x_2} \left(-{a_1} {d_1} {x_1}^2-{a_1} {d_1}{x_2}^2-{a_1} {d_1} {x_3}^2-{a_2} {d_1} {x_1}^2+{a_2}
{d_1} {x_2}^2-{a_2} {d_1} {x_3}^2-{a_2} {d_2}{x_1}^2+{a_2} {d_2} {x_2}^2 \right. \nonumber\\
&& \left. -{a_2} {d_2} {x_3}^2-{a_2}{d_3} {x_1}^2+{a_2} {d_3} {x_2}^2-{a_2} {d_3}{x_3}^2+{d_1} {x_1} {x_3}+{d_1} {x_2} {x_3}+2 {d_2}{x_1} {x_3}+{d_3} {x_1} {x_2}+{d_3} {x_1} {x_3}\right) \nonumber\\
\dot{x}_3 &=& {x_3} \left(-{a_1} {d_1} {x_1}^2-{a_1} {d_1}{x_2}^2-{a_1} {d_1} {x_3}^2-{a_3} {d_1} {x_1}^2-{a_3}
{d_1} {x_2}^2+{a_3} {d_1} {x_3}^2-{a_3} {d_2}{x_1}^2-{a_3} {d_2} {x_2}^2 \right. \nonumber\\
&& \left. +{a_3} {d_2} {x_3}^2-{a_3}{d_3} {x_1}^2-{a_3} {d_3} {x_2}^2+{a_3} {d_3}{x_3}^2+{d_1} {x_1} {x_2}+{d_1} {x_2} {x_3}+{d_2}{x_1} {x_2}+{d_2} {x_1} {x_3}+2 {d_3} {x_1}{x_2}\right)
\end{eqnarray}
If we denote by $P_i(x_1, x_2, x_3)$ the $\dot{x}_i$, for $i=1,2,3$ respectively, then the degree of the vector field $X = \sum_{i=1}^3P_i(x_1, x_2, x_3)\frac{\partial}{\partial x_i}$ is $3$.  Now 
we study the system (\ref{systemWa}) at infinity.  We apply the Poincar\'e compactification to the above system written in the chart $(U_1, \varphi_1)$ as follows:
\begin{eqnarray}\label{systemWallach}
\dot{z}_1 &=& {z_1} ({d_1}+{d_2}+{d_3}) \left(\left({z_1}^2-1\right)({a_1}+{a_2})+{z_2}^2 ({a_1}-{a_2})-{z_1} {z_2}+{z_2}\right) \nonumber\\
\dot{z}_2 &=& {z_2} ({d_1}+{d_2}+{d_3}) \left({a_1}\left({z_1}^2+{z_2}^2-1\right)+{a_3}
   \left(-{z_1}^2+{z_2}^2-1\right)+{z_1} (-{z_2})+{z_1}\right)  \nonumber\\
\dot{z}_3 &=& {z_3} \left({a_1} \left(2 {d_1}\left({z_1}^2+{z_2}^2\right)+{d_2}\left({z_1}^2+{z_2}^2-1\right)+{d_3}\left({z_1}^2+{z_2}^2-1\right)\right) \right. \nonumber \\
&& \left. -{z_2} ({z_1} (2{d_1}+{d_2}+{d_3})+{d_2})-{d_3} {z_1}\right).
\end{eqnarray}

In order to find the singularities at infinity of the above system we set $z_3 = 0$.  We will substitute into (\ref{systemWallach}), the values of the dimensions $d_1, d_2, d_3$ and $a_1, a_2, a_3$ from Tables 1 and 2 respectively.  Then we have the following:\\
For the generalized Wallach space {\small\small\bf GWS.1}\ $\SO(k+l+m)/\SO(k)\times\SO(l)\times\SO(m)$  system (\ref{systemWallach}) comes 
\begin{eqnarray}\label{WallachSO}
\dot{z}_1 &=& {z_1} (-k (l+m)-l m) \left(2 ({z_1}-1) {z_2}(k+l+m-2)+\left({z_1}^2-1\right) (-l-m)+{z_2}^2 (l-m)\right)\nonumber \\
\dot{z}_2 &=& {z_2} (-k (l+m)-l m) \left(k \left(2 {z_1}{z_2}+({z_1}-1)^2-{z_2}^2\right)+2 (l-2) {z_1} ({z_2}-1)\right. \nonumber\\
&&\left. -m\left(-2 {z_1} {z_2}+{z_1} ({z_1}+2)+{z_2}^2-1\right)\right).
\end{eqnarray}

Next, we will study the case $l = m$.

\begin{prop}\label{Propo1}
For the generalized Wallach space $\SO(k+2m)/\SO(k)\times\SO(m)\times\SO(m)$, the normalized Ricci flow has for  $m>2$ and $k^2 > 4(m-1)$,  precisely four singularities at infinity.
\end{prop}
\begin{proof}
We have $m=l$ so  system (\ref{systemWallach}) can be written as
\begin{eqnarray}\label{WallachSOkl}
\dot{z}_1 = f_1(z_1, z_2) &=& 2 m ({z_1}-1) {z_1} (2 k+m) (m ({z_1}+1)-{z_2} (k+2 m-2))=0 \nonumber\\
\dot{z}_2 = f_2(z_1, z_2) &=& -m {z_2} (2 k+m) (k (2 {z_1}{z_2}+({z_1}-1)^2-{z_2}^2)+4 (m-1) {z_1} {z_2} \nonumber\\
&& -m {z_1}({z_1}+4)-m {z_2}^2+m+4 {z_1})=0.
\end{eqnarray}

From the first equation we obtain $z_1=1$ and  by substituting into the (\ref{WallachSOkl}) we 
obtain the solutions
$$z_2 = \frac{k+2 m-2\pm \sqrt{k^2-4 m+4}}{k+m},
$$  
hence we obtain two solutions $(1, z_2)$.  Now, if $ (m ({z_1}+1)-{z_2} (k+2 m-2))=0$ then it is possible to obtain explicit solutions for $z_1, z_2$, but these are complicated expressions involving radicals.  So we use Gr\"obner bases to prove existence of at least two solutions.  We consider a polynomial ring $R = \bb{Q}[z, z_1, z_2]$ and an ideal $I$ generated by the polynomials $\{f_1(z_1, z_2),$ $f_2(z_1, z_2),$  $z\, z_1\, z_2\, (z_1-1) - 1\}$.  We take a lexicographic ordering $>$ with $z > z_1 > z_2$ for a monomial ordering on $R$.  Then, with the aid of computer, we see that a Gr\"obner basis for the ideal $I$ contains the polynomial $H(z_2)$ of $z_2$, which is given by
{\small \begin{eqnarray*} 
&& H(z_2) = 8 k^2 m^3+ (-8 k^3 m^2-28 k^2 m^3+24 k^2 m^2-28 k m^4+44 k m^3-16 k m^2-8m^5+16 m^4 -8 m^3){z_2}\\
&&  +(2 k^4 m+11 k^3 m^2-8 k^3 m+19 k^2 m^3-28 k^2m^2+8 k^2 m+13 k m^4-28 k m^3+12 k m^2+3 m^5\\
&& -8 m^4 +4 m^3){z_2}^2 +12 k m^4+8 k m^3+4m^5-4 m^4.
\end{eqnarray*} }
We compute $H(0) = 8 k^2 m^3+12 k m^4-8 k m^3+4 m^5-4 m^4 >0$, $H(1) = 2 k^4 m+3 k^3 m^2-8 k^3 m-k^2 m^3-4 k^2 m^2+8 k^2 m-3 k m^4+8 k m^3-4 k m^2-m^5+4 m^4-4 m^3 < 0$ and $H(2) =8 k^4 m+28 k^3 m^2-32 k^3 m+28 k^2 m^3-64 k^2 m^2+32 k^2 m+8 k m^4-32 k m^3+16 k m^2-4 m^4 $ which is positive for $m > 2$.  Thus, there exist at least two solutions $z_2 = \al$ and $z_2 = \beta$ of $H(z_2) = 0$ with $0 < \al < 1$ and $1 < \beta < 2$.  Now, for positive values $z_2 = \al, \beta$ we obtain from $ (m ({z_1}+1)-{z_2} (k+2 m-2))=0$ real values $z_1 = \gamma, \delta$ as solutions of (\ref{WallachSOkl}).  Next, we consider a lexicographic ordering $>$ with $z > z_2 > z_1$ for a monomial ordering on the ring $R$. Then, a Gr\"obner basis for the ideal $I$ contains the polynomial $G(z_1)$ of $z_1$, which can be written as
{\small \begin{eqnarray*}  
&& G(z_1) = (6 k^4+75 k^3+(19 k^2+128 k+178) (m-3)^3+285 k^2+(11 k^3+143 k^2+462k+414) (m-3)^2\\
&& +(2 k^4+58 k^3+353 k^2+720 k+459) (m-3)+(13 k+37)(m-3)^4+405 k+3 (m-3)^5+189)z_1^2 + (-12 k^4\\
&& -126 k^3+(-46 k^2-352 k-556) (m-3)^3-558 k^2+(-22 k^3-322k^2-1220 k-1316) (m-3)^2+(-4 k^4-108 k^3\\
&& -738 k^2-1888 k-1578)(m-3)+(-38 k-118) (m-3)^4-1110 k-10 (m-3)^5-774  )z_1 + 6 k^4+75 k^3+(19 k^2\\
&& +128 k+178) (m-3)^3+285 k^2+(11 k^3+143 k^2+462k+414) (m-3)^2+(2 k^4+58 k^3+353 k^2+720 k\\
&& +459) (m-3)+(13 k+37)(m-3)^4+405 k+3 (m-3)^5+189.
\end{eqnarray*} }
Then we see that, for $m\geq 3$, the coefficients of the polynomial $G(z_1)$ are positive for even degree terms and negative for odd degree terms.  Thus, if the equation $G(z_1) = 0$ has real solutions, then these are all positive.  In particular, $z_1 = \gamma, \delta$ are positive.
So far we have proved existence of at least four solutions.  However, by Theorem \ref{general} the number of invariant Einstein metrics is at most four, hence the result follows.
\end{proof}

\begin{remark} 
{\rm If $k^2=4(m-1)$, then $k$ and $m$ are of the form $(k,m)=(2\rho, \rho^2+1)$, $\rho >1$, and according to Example 1, p. 56 in \cite{LoNiFi}, there exist precisely three singularities at infinity.
Also, if $m=2$ then $k=2$, this corresponds to the space $\SO(6)/\SO(2)\times\SO(2)\times\SO(2)$.  It easy to see that the system (\ref{WallachSOkl})  has a unique singularity at infinity.  This singularity corresponds to the unique invariant Einstein metric (cf. \cite{LoNiFi}, \cite{CN}). }
\end{remark}

\noindent For the space {\small\small\bf GWS.2}\ $\SU(k+l+m)/\SU(k)\times\SU(l)\times\SU(m)$ system (\ref{systemWallach}) becomes:
\begin{eqnarray*}
&&\dot{z}_1 = {z_1} (-k (l+m)-l m) \left(2 ({z_1}-1) {z_2}
   (k+l+m)+\left({z_1}^2-1\right) (-l-m)+{z_2}^2 (l-m)\right),\\
&&\dot{z}_2 = {z_2} (-k (l+m)-l m) (2 {z_1} {z_2} (k+l+m)-{z_2}^2 (k+m)+k
   {z_1}^2-2 k {z_1}+k-2 l {z_1}-m {z_1}^2-2 m {z_1}+m).  
\end{eqnarray*}
The solutions are: $((k + m)/(k + l), (l + m)/(k + l)), ((k + m)/(k + l), (2 k + l + m)/(k + l)), ((k + 2 l + m)/(k + l), (l + m)/(k + l))$ and $((k + m)/(k + l + 2 m), (l + m)/(k + l + 2 m))$.

\smallskip\smallskip
\noindent For the generalized Wallach space {\small\small\bf GWS.3}\ $\Sp(k+l+m)/\Sp(k)\times\Sp(l)\times\Sp(m)$ 
the system (\ref{systemWallach}) comes: 
\begin{eqnarray*}
\dot{z}_1 &=& {z_1} (-k (l+m)-l m) \left(2 ({z_1}-1) {z_2}(k+l+m+1)+\left({z_1}^2-1\right) (-l-m)+{z_2}^2 (l-m)\right)\\
\dot{z}_2 &=& {z_2} (-k (l+m)-l m) \left(k \left(2 {z_1}{z_2}+({z_1}-1)^2-{z_2}^2\right)+2 (l+1) {z_1} ({z_2}-1)\right.\\
&&\left. -m\left(-2 {z_1} {z_2}+{z_1} ({z_1}+2)+{z_2}^2-1\right)\right).
\end{eqnarray*}
By similar method as in Proposition \ref{Propo1} we can prove the following:
\begin{prop}
For the generalized Wallach space $\Sp(k+2m)/\Sp(k)\times\Sp(m)\times\Sp(m)$, the normalized Ricci flow has precisely four singularities at infinity. 
\end{prop}

It is easy to see that for some values of $k,l$ and $m$ the above system has four solutions.  For example,
\begin{itemize}
\item[(1)] For $(k, l, m) = (1, 2, 3)$ the solutions are: $(3.26361, 1.60389), (1.30670, 3.18223), (1.23251, 1.39606)$ and $(0.38050, 0.46780)$.
\item[(2)] For $(k, l, m) = (2, 5, 7)$ the solutions are: $(2.94748, 1.67504), (1.27217, 2.71689), (1.24716, 1.53155)$ and $(0.40168, 0.52944)$.
\end{itemize}



\noindent{\small\small\bf GWS.4}\ $\SU(2l)/\U(l)$
\begin{eqnarray}\label{WallSU}
\dot{z}_1 &=& 2 \left(3 l^2-1\right) {z_1} \left(l ({z_1}-1) ({z_1}-2{z_2}+1)+{z_2}^2\right)  \nonumber\\
\dot{z}_2 &=& \left(3 l^2-1\right) {z_2} \left(2 l \left(-2 {z_1}
   ({z_2}-1)+{z_2}^2-1\right)+{z_1}^2+{z_2}^2-1\right).
\end{eqnarray}

\begin{prop}
For the generalized Wallach space $\SU(2l)/\U(l)$, the normalized Ricci flow has for $l\geq 2$ at least two singularities at infinity.
\end{prop}
\begin{proof}
We consider a polynomial ring $R = \bb{Q}[z, z_1, z_2]$ and an ideal $I$ generated by the polynomials $\{\dot{z}_1 = \left(3 l^2-1\right) {z_1} \left(l ({z_1}-1) ({z_1}-2{z_2}+1)+{z_2}^2\right)$ $\dot{z}_2=\left(3 l^2-1\right) {z_2} \left(2 l \left(-2 {z_1}
   ({z_2}-1)+{z_2}^2-1\right)+{z_1}^2+{z_2}^2-1\right) ,$  $z\, z_1\, z_2-1\}$.  We take a lexicographic ordering $>$ with $z_2 > z_1$ for a monomial ordering on $R$.  Then, with the aid of computer we see that a Gr\"obner basis for the ideal $I$ contains the polynomial $F(z_1)$ of $z_1$, which is given by
\begin{eqnarray*}
F(z_1) &=& (3l^2-1) 12 l^4 {z_1}^4-48 l^4 {z_1}^3+72 l^4 {z_1}^2-48 l^4 {z_1}+12 l^4-20 l^3{z_1}^4+40 l^3 {z_1}^3-40 l^3 {z_1}+20 l^3\\
&& +7 l^2 {z_1}^4+8 l^2{z_1}^3-34 l^2 {z_1}^2+8 l^2 {z_1}+7 l^2+2 l {z_1}^4-8 l {z_1}^3+8l {z_1}-2 l-{z_1}^4+2 {z_1}^2-1.
\end{eqnarray*}
We compute $F(0) = 12 l^4+20 l^3+7 l^2-2 l-1 >0$, $H(1) = -4 l^2<$ and $H(4) = 972 l^4-2700 l^3+1799 l^2+30 l-225$, which is positive for $l\geq 2$.  Thus, there exist at least two solutions $z_1 = \al$ and $z_1 = \beta$ of $F(z_1) = 0$ with $0 < \al_1 < 1$ and $1 < \beta_1 < 4$.  The Gr\"obner basis for the ideal $I$ contains also the polynomial
$
 (3l^2-1)z_2 -  (3l^2-1)a(z_1, k) = 0.
$ 
Thus for the positive values $z_1=\al_1, \beta_1$ found above, we obtain real values $z_2 = \al_2, \beta_2$ as solutions of the system (\ref{WallSU}).  Now we consider a lexicographic ordering $>$ with $z_1 > z_2$ for a monomial ordering on the ring $R$. Then, a Gr\"obner basis for the ideal $I$ contains the polynomial polynomial $(3l^2-1) G(z_2)$, where $G(z_2)$ is a polynomial of of $z_2$, which can be written as
\begin{eqnarray*}
G(z_2) &=& (12 (l-2)^4+76 (l-2)^3+175 (l-2)^2+174 (l-2)+63)z_2^4 \\
&& (-48 (l-2)^4-336 (l-2)^3-868 (l-2)^2-980 (l-2)-408)z_2^3 \\
&& (72 (l-2)^4+540 (l-2)^3+1508 (l-2)^2+1856 (l-2)+848)z_2^2 \\
&& (-48 (l-2)^4-376 (l-2)^3-1104 (l-2)^2-1440 (l-2)-704)z_2 \\
&& 12 (l-2)^4+96 (l-2)^3+288 (l-2)^2+384 (l-2)+192.
\end{eqnarray*}
Then we see that, for$l\geq 2$,the coefficients of the polynomial $G(z_2)$ are positive for even degree terms and negative for odd degree terms.  Thus if the equation $G(z_2) = 0$ has real solutions, then these are all positive.  In particular, $z_2 = \al_2, \beta_2$ are positive.
\end{proof}


\noindent{\small\small\bf GWS.5}\ $\SO(2l)/\U(1)\times\U(l-1)$
\begin{eqnarray*}
\dot{z}_1 &=& 2 (l+2) ({z_1}-1) {z_1} ((l-2) ({z_1}+1)-2 (l-1) {z_2}), \\
\dot{z}_2 &=& (l+2) {z_2} \left(l \left(-4 {z_1} {z_2}+{z_1}
   ({z_1}+4)+{z_2}^2-1\right)-4 {z_1} ({z_1}-{z_2}+1)\right)
\end{eqnarray*}
The solutions are: $(1, 2), (1, (2l-4)/l), (l/(3l-4), 2(l-2)/(3l-4))$ and $((3l-4)/l, 2(l-2)/l)$.

\noindent{\small\small\bf GWS.6}\ $\E_6/\SU(4)\times\SU(2)\times\SU(2)\times\U(1)$
\begin{eqnarray*}
\dot{z}_1 = 28 ({z_1}-1) {z_1} ({z_1}-2 {z_2}+1), \ \ \dot{z}_2 = (14/3){z_2} \left({z_1}^2-12 {z_1} ({z_2}-1)+5 {z_2}^2-5\right).
\end{eqnarray*}  
The solutions are: $(0.6, 0.8)$ and $(1.66667, 1.33333)$.

\noindent{\small\small\bf GWS.7}\ $\E_6/\SO(8)\times\U(1)\times\U(1)$
\begin{eqnarray*}
\dot{z}_1 = 16 ({z_1}-1) {z_1} ({z_1}-3 {z_2}+1), \ \ \dot{z}_2 = 16 {z_2} \left(-3 {z_1} ({z_2}-1)+{z_2}^2-1\right).
\end{eqnarray*}
The solutions are: $(1, 1)$, $(2, 1)$, $(1, 2)$ and $(1/2, 1/2)$.

\noindent{\small\small\bf GWS.8}\ $\E_6/\Sp(3)\times\SU(2)$
\begin{eqnarray*}
\dot{z}_1 = (27/4){z_1} \left(3 {z_1}^2-8 {z_1} {z_2}+{z_2} ({z_2}+8)-3\right), \ \ \dot{z}_2 = -(9/4){z_2} \left({z_1}^2+24 {z_1} ({z_2}-1)-13 {z_2}^2+13\right).
\end{eqnarray*}
The solutions are: $(0.864003, 0.483834)$ and $(1.46177, 1.884488)$.

\noindent{\small\small\bf GWS.9}\ $\E_7/\SO(8)\times\SU(2)\times\SU(2)\times\SU(2)$
\begin{eqnarray*}
\dot{z}_1 = (32/3)({z_1}-1) {z_1} (4 {z_1}-9 {z_2}+4), \ \  \dot{z}_2 = (32/3)({z_2}-1) {z_2} (-9 {z_1}+4 {z_2}+4).
\end{eqnarray*}
The solutions are: $(1, 1)$, $(1.25, 1)$, $(1, 1.25)$ and $(0.8, 0.8)$.


\noindent{\small\small\bf GWS.11}\ $\E_7/\SO(8)$
\begin{eqnarray*}
\dot{z}_1 = (35/3)({z_1}-1) {z_1} (5 {z_1}-9 {z_2}+5), \ \  \dot{z}_2 = (35/3)({z_2}-1) {z_2} (-9 {z_1}+5 {z_2}+5).
\end{eqnarray*}
The solutions are: $(1, 1)$, $(1, 0.8)$, $(0.8, 1)$ and $(1.25, 1.25)$.

\noindent{\small\small\bf GWS.12}\ $\E_8/\SO(12)\times\SU(2)\times\SU(2)$
\begin{eqnarray*}
\dot{z}_1 = (176/5)({z_1}-1) {z_1} (2 {z_1}-5 {z_2}+2), \ \ \dot{z}_2 = -(176/5){z_2} \left({z_1}^2+15 {z_1} ({z_2}-1)-7 {z_2}^2+7\right).
\end{eqnarray*}
The solutions are: $(1, 1.456083)$ and $(1, 0.686773)$.

\noindent{\small\small\bf GWS.13}\ $\E_8/\SO(8)\times\SO(8)$
\begin{eqnarray*}
\dot{z}_1 = (64/5)({z_1}-1) {z_1} (8 {z_1}-15 {z_2}+8), \ \  \dot{z}_2 = (64/5)({z_2}-1) {z_2} (-15 {z_1}+8 {z_2}+8).
\end{eqnarray*}
The solutions are: $(1, 1)$, $(1, 0.875)$, $(0.875, 1)$ and $(1.142857, 1.142857)$.

\noindent{\small\small\bf GWS.14}\ $\F_4/\SO(5)\times\SU(2)\times\SU(2)$
\begin{eqnarray*}
\dot{z}_1 = 4 ({z_1}-1) {z_1} (5 {z_1}-9 {z_2}+5), \ \  \dot{z}_2 = 2 {z_2} \left(-18 {z_1} {z_2}+3 {z_1} ({z_1}+6)+7{z_2}^2-7\right)
\end{eqnarray*}
The solutions are: $(0.485288, 0.825160)$ and $(2.060629, 1.700349)$.

\noindent{\small\small\bf GWS.15}\ $\F_4/\SO(8)$
\begin{eqnarray*}
\dot{z}_1 = (8/3)({z_1}-1) {z_1} (2 {z_1}-9 {z_2}+2), \ \  \dot{z}_2 = (8/3)({z_2}-1) {z_2} (-9 {z_1}+2 {z_2}+2).
\end{eqnarray*}
The solutions are: $(1, 1)$, $(3.5, 1)$, $(1, 3.5)$ and $(0.285714, 0.285714)$.

Since we work on the chart $(U_1, \varphi_1)$ that corresponds to the plane $y_1 = 1$, we consider metrics which are of the form $(1, \alpha, \beta)$, where $\alpha, \beta$ are the solutions of system (\ref{systemWallach}) for $z_3 = 0$. 
These metrics are invariant Einstein metrics on $M$.  Thus, we have proved the following:

\begin{theorem}
Let $G/H$ be a generalized Wallach space.  The normalized Ricci flow, on the space of invariant Riemannian metrics on $G/H$, possesses a finite number of singularities 
at infinity.  These fixed points correspond to the $G$-invariant Einstein metrics on $G/H$ (cf. Subsection 3.1). 
\end{theorem}

\subsection{Ricci flow for the Stiefel manifold $V_2\bb{R}^n$}
By using the Ricci components (\ref{ricciV2}) and scalar curvature (\ref{scalV2}) of the metric (\ref{metricSt}) it is easy to see that, for the Stiefel manifold $V_2\bb{R}^n$ system (\ref{normalRicci}) reduces to 
\begin{eqnarray*}
\dot{x}_{0} &=& \frac{(n-2) {x_0}^2+(n-2) {x_0} ({x_1}+{x_2})-(n-1)
   ({x_1}-{x_2})^2}{(2 n-3) {x_1} {x_2}}\\
\dot{x}_1 &=&  \frac{(5-3 n) {x_0}^2+2 (n-2) {x_0} ((n-2) {x_1}+(3 n-5)
   {x_2})+({x_1}-{x_2}) ((n-1) {x_1}+(3 n-5) {x_2})}{2 (n-2) (2 n-3)
   {x_0} {x_2}}\\
\dot{x}_2 &=& \frac{(5-3 n) {x_0}^2+2 (n-2) {x_0} ((3 n-5) {x_1}+(n-2)
   {x_2})-({x_1}-{x_2}) ((3 n-5) {x_1}+(n-1) {x_2})}{2 (n-2) (2 n-3)
   {x_0} {x_1}}.
\end{eqnarray*}
We observe that the above system is not a polynomial system, therefore we cannot apply the Poincar\'e compactification directly.  We multiply it by the factor $2(n - 2)(2n - 3)x_0x_1x_2$, with this multiplication  will only change the time of parametrization of the orbits and not the structure of the phase portrait.  After that we take the following system:
\begin{eqnarray}\label{systemSt}
\dot{x}_0 &=& (n-2) \left(n {x_0}^2+n {x_0} {x_1}+n {x_0} {x_2}-n{x_1}^2+2 n {x_1} {x_2}\right. \nonumber\\
&&\left. -n {x_2}^2-2 {x_0}^2-2 {x_0}
   {x_1}-2 {x_0} {x_2}+{x_1}^2-2 {x_1} {x_2}+{x_2}^2\right){x_0} \nonumber\\
\dot{x}_1 &=& 1/2\left(2 n^2 {x_0} {x_1}+6 n^2 {x_0} {x_2}-3 n {x_0}^2-8 n{x_0} {x_1}-22 n {x_0} {x_2}\right.\nonumber \\
&& \left. +n {x_1}^2+2 n {x_1} {x_2}-3 n{x_2}^2+5 {x_0}^2+8 {x_0} {x_1}+20 {x_0} {x_2}-{x_1}^2-4{x_1} {x_2}+5 {x_2}^2\right) {x_1} \nonumber\\
\dot{x}_2 &=& 1/2  \left(6 n^2 {x_0} {x_1}+2 n^2 {x_0} {x_2}-3 n {x_0}^2-22 n{x_0} {x_1}-8 n {x_0} {x_2}\right.\nonumber\\
&& \left. -3 n {x_1}^2+2 n {x_1} {x_2}+n{x_2}^2+5 {x_0}^2+20 {x_0} {x_1}+8 {x_0} {x_2}+5 {x_1}^2-4{x_1} {x_2}-{x_2}^2\right){x_2}.
\end{eqnarray}

If we denote by $P_i(x_0, x_1, x_2)$ the $\dot{x}_i$, for $i=0,1,2$ respectively, then the degree of the vector field $X = \sum_{i=0}^2P_i(x_0, x_1, x_2)\frac{\partial}{\partial x_i}$ is $3$.
From the above system it is easy to see the following lemma:

\begin{lemma}
The coordinate planes along with the straight line $\gamma(t) = \left((2(n-2)/(n-1))t, t, t\right)$, remain invariant under the flow defined by the system of (\ref{systemSt}).
\end{lemma}
   
To study the singularities at infinity of (\ref{systemSt}), we should write this system in the local charts of the Poincar\'e compactification.  Because we are interesting for positive values of $x_0, x_1$ and $x_2$ we study the behavior of the previous system only on the chart $(U_1, \varphi_1)$.  Therefore from (\ref{Xartis1}) we have:
\begin{eqnarray}\label{RiciSt1}
\dot{z}_1 &=& 1/2 (2 n-3) {z_1} \left((n-1) \left({z_1}^2-1\right)-2 (n-2) ({z_1}-1)
   {z_2}+(n-3) {z_2}^2\right) \nonumber\\
\dot{z}_2 &=& 1/2 (2 n-3) {z_2} \left(n \left(-2 {z_1} {z_2}+{z_1}({z_1}+2)+{z_2}^2-1\right)-3 {z_1}^2+4 {z_1}
({z_2}-1)-{z_2}^2+1\right) \nonumber\\
\dot{z}_3 &=& (n-2) {z_3} \left(n \left({z_1}^2-{z_1} (2 {z_2}+1)+({z_2}-1)
   {z_2}-1\right)-({z_1}-{z_2})^2+2 ({z_1}+{z_2}+1)\right).
\end{eqnarray} 
In order to find the singularities at the infinity of the above system we set $z_3 = 0$.  Then it is easy to see that the system $\{\dot{z}_1 = 0, \dot{z}_2 = 0\}$ has only one solution namely $(z_1, z_2) = ({(n-1)}/{2(n-2)}, {(n-1)}/{2(n-2)})$.  

Since we work on the chart $(U_1, \varphi_1)$ that corresponds to the plane $y_1 = 1$, we consider metrics whose are the form $(1, \alpha, \beta)$, where $\alpha, \beta$ are the solution of the system (\ref{RiciSt1}) for $z_3 = 0$.  These metrics are invariant Einstein metrics on $M$.  We have thus proved the following.

\begin{theorem}
Let $G/H$ be the Stiefel manifold $V_2\bb{R}^n$.  The normalized Ricci flow on the space of invariant Riemannian metrics on $G/H$, possesses exactly one singularity at infinity.  This fixed point corresponds to the unique (up to scale) $G$-invariant Einstein metric on $G/H$ (cf. Theorem 4.3). 
\end{theorem}

\subsection{Ricci flow on the Stiefel manifold $V_{1+k_2}\bb{R}^n$}
We study the behavior of the normalized Ricci flow for the Stiefel manifold $V_{1+k_2}\bb{R}^n$.  For this case  
we take the Ricci components (\ref{eq19}) and scalar curvature (\ref{scalarSt2}) of the metric (\ref{metric2}).  Then, system (\ref{normalRicci}) reduces to 
{ \begin{eqnarray*}
&&\dot{x}_{2} =
-\frac{1}{{2({k_2}+1)({k_2}-2 n+2)(n-2) {x_{12}}^2 {x_{13}}{x_{23}}^2 }}  
 \Big\lbrace  4 ({x_{23}}^2 (-{x_{13}} {x_2}^2({k_2}^2-({k_2}+1) n+{k_2}+1)  \\
&&  +{k_2} {x_{12}}^2 {x_{13}} ({k_2} (n-3)+n)+2 (n-2) {x_{12}} {x_2} ({x_12}(-{k_2}+n-1)+{k_2} {x_{13}}))+{x_{12}}^2 {x_{13}} {x_2}^2({k_2}\\
&& -n+1) ({k_2}^2-({k_2}+1) n+{k_2}+1)+{k_2}{x_{12}} {x_2} {x_{23}} ({k_2}-n+1) (-2 (n-2) {x_{12}}{x_{13}}+{x_{12}}^2+{x_{13}}^2)\\
&& +{k_2} {x_{12}} {x_2}{x_{23}}^3 ({k_2}-n+1))\Big\rbrace.
\\
&&\dot{x}_{12} = -\frac{1}{{2({k_2}+1) (n-2)({k_2}-2 n+2) {x_{12}} {x_{13}} {x_2} {x_{23}}^2 }}\Big\lbrace {x_{23}} ({x_{12}}^3 (-{k_2}+n-1)+{x_{12}} ({x_{13}}^2({k_2}-n+1)\\
&& + {x_{23}}^2 ({k_2}-n+1)+2 (n-2){x_{13}}{x_{23}})-({k_2}-1) {x_{13}} {x_2} {x_{23}})+2 {k_2}{x_{12}} {x_2} {x_{23}} ({k_2}-n +1) (-2 (n\\
&& -2) {x_{12}}{x_{13}}+{x_{12}}^2+{x_{13}}^2)+{x_{23}}^2 (4 (n-2) {x_{12}}{x_2} ({x_{12}} (-{k_2}+n-1)+{k_2} {x_{13}})+({k_2}-2)({k_2}-1) {k_2} \\
&& {x_{12}}^2 {x_{13}}-({k_2}-1) {k_2} {x_{13}}{x_2}^2)+({k_2}-1) {k_2} {x_{12}}^2 {x_{13}} {x_2}^2
({k_2}-n+1)+2 {k_2}{x_{12}} {x_2} {x_{23}}^3 ({k_2}-n+1) \Big\rbrace.
\\
&&\dot{x}_{13} = -\frac{1}{2({k_2}+1) (n-2)({k_2}-2 n+2) {x_{12}}^2 {x_2} {x_{23}}^2 } \Big\lbrace  2 {k_2} {x_{12}} {x_2} {x_{23}} ({k_2}-n+1) (-2 (n-2){x_{12}} {x_{13}}+{x_{12}}^2\\
&& +{x_{13}}^2)+{x_{12}} {x_{23}}({k_2} ({x_{12}}^2-{x_{13}}^2+{x_{23}}^2)-2 (n-2){x_{12}} {x_{23}})+{x_{23}}^2 (4 (n-2) {x_{12}} {x_2}({x_{12}} (-{k_2}+n\\
&& -1)+{k_2} {x_{13}})+({k_2}-2) ({k_2}-1){k_2} {x_{12}}^2 {x_{13}}-({k_2}-1) {k_2} {x_{13}}{x_2}^2)+({k_2}-1) {k_2} {x_{12}}^2 {x_{13}} {x_2}^2 ({k_2}-n\\
&& +1)+2 {k_2} {x_{12}} x_{2} {x_{23}}^3 ({k_2}-n+1)  \Big\rbrace.
\\
&&\dot{x}_{23} = -\frac{1}{2 ({k_2}+1) (n-2)({k_2}-2 n+2) {x_{12}}^2 {x_{13}} {x_2} {x_{23}}} \Big\lbrace   {x_{12}} {x_2} {x_{23}}^3 ({k_2}^2-{k_2}+2n-2)+{x_{12}} {x_2} {x_{23}} ({k_2} (3 {k_2}\\
&& -4 n+5)-2 n+2)(-2 (n-2) {x_{12}} {x_{13}}+{x_{12}}^2+{x_{13}}^2)+{x_{23}}^2(4 (n-2) {x_{12}} {x_2} ({x_{12}} (-{k_2}+n-1)\\
&&+{k_2}{x_{13}})+({k_2}-2) ({k_2}-1) {k_2} {x_{12}}^2{x_{13}}-({k_2}-1){k_2} {x_{13}} {x_2}^2)\\
&& +({k_2}-1){x_{12}}^2 {x_{13}} {x_2}^2 (2 ({k_2}+1)^2-(3 {k_2}+2)n) \Big\rbrace.  
\end{eqnarray*}  }
The above system is not polynomial, hence in order to apply the Poincar\'e compactification, we multiply it by the factor $2({k_2}+1)(k_2+k_3-1)(-{k_2}-2k_3){x_{12}}^2 {x_{13}} {x_2} {x_{23}^2}$.  Then we obtain the following system:
{ \begin{eqnarray}\label{systemStiefel}
\dot{x}_2 &=& 2 {x_2} ({k_2} {k_3} {x_{12}} {x_2} {x_{23}} (2 {x_{12}}{x_{13}} ({k_2}+{k_3}-1)-{x_{12}}^2-{x_{13}}^2) +{x_{23}}^2(({k_2}-2) {x_{12}}^2 {x_{13}} ({k_2}({k_2}+{k_3}) \nonumber \\
&& +{k_3})+2 {x_{12}} {x_2} ({k_2}+{k_3}-1)({k_2} {x_{13}}+{k_3} {x_{12}})+{x_{13}} {x_2}^2 ({k_2}{k_3}+{k_2}+{k_3}))+ \nonumber\\
&& {k_3} {x_{12}}^2 {x_{13}} {x_2}^2({k_2} {k_3}+{k_2}+{k_3})-{k_2} {k_3} {x_{12}} {x_2}{x_{23}}^3) \nonumber\\
\dot{x}_{12} &=& {x_{12}} ({k_2}^3 {x_{13}} {x_{23}}^2 ({x_{12}}^2+2 {x_{12}}{x_2}-{x_2}^2)+{k_2}^2({k_3} {x_{12}} {x_2} {x_{23}}({x_{12}}^2+4 {x_{12}} {x_{13}}-{x_{13}}^2) \nonumber \\
&& -{k_3}{x_{12}}^2 {x_{13}} {x_2}^2-{x_{13}} {x_{23}}^2 (3 {x_{12}}-{x_2})({x_{12}}-(2 {k_3}+1) {x_2})-{k_3} {x_{12}} {x_2}{x_{23}}^3) \nonumber \\
&& +{k_2} (2 {k_3}^2 {x_{12}} {x_2} {x_{23}} ({x_{12}}^2+2 {x_{12}} {x_{13}}-({x_{13}}-{x_{23}})^2)+{k_3} {x_{12}} {x_2} ({x_{12}} {x_{13}} {x_2}+2 {x_{23}}^2 (2 {x_{12}}+3 {x_{13}}) \nonumber\\
&& -{x_{23}} ({x_{12}}+{x_{13}}) ({x_{12}}+3 {x_{13}})-3 {x_{23}}^3)+2 {x_{13}} {x_{23}}^2 ({x_{12}}^2-3{x_{12}} {x_2}+{x_2}^2))+2 {k_3} {x_2} {x_{23}}({k_3}  \nonumber \\
&& {x_{12}} ({x_{12}}^2+2 {x_{12}}{x_{23}}-({x_{13}}-{x_{23}})^2)-2 {x_{12}} {x_{23}}({x_{12}}+{x_{13}})+{x_{13}} {x_2} {x_{23}})) \nonumber\\
\dot{x}_{13} &=& {x_{13}} ({k_2}^3 {x_{12}} {x_{23}} ({x_{12}}^2 (-{x_2})+{x_{12}} {x_{23}} ({x_{13}}+2 {x_2})+{x_2}({x_{13}}-{x_{23}}) ({x_{13}}+{x_{23}}))-{k_2}^2 ({k_3}{x_{12}} {x_2} (2 {x_{23}} \nonumber\\
&& ({x_{12}}^2-2 {x_{12}}{x_{13}}-{x_{13}}^2)+{x_{12}} {x_{13}} {x_2}-6 {x_{12}}{x_{23}}^2+2 {x_{23}}^3)+{x_{23}} ({x_{13}} {x_{23}}({x_{12}}-{x_2}) (3 {x_{12}}-{x_2})\nonumber\\
&& +{x_{12}} {x_2} ({x_{12}}-{x_{13}}) ({x_{12}}+{x_{13}})+{x_{12}} {x_2}{x_{23}}^2))+{k_2} (4 {k_3}^2 {x_{12}}^2 {x_2} {x_{23}} ({x_{13}}+{x_{23}})+{k_3} {x_{12}} {x_2} ({x_{12}} {x_{13}} {x_2} \nonumber\\
&& +2 {x_{23}}^2 (3 {x_{12}}+2 {x_{13}})-4 {x_{12}}{x_{23}} ({x_{12}}+{x_{13}})-4 {x_{23}}^3)+{x_{23}}^2 (2{x_{12}}^2 {x_{13}}-2 {x_{12}} {x_2} ({x_{12}}+2 {x_{13}})\nonumber\\
&& +{x_{13}}{x_2}^2))+8 ({k_3}-1) {k_3} {x_{12}}^2 {x_2}{x_{23}}^2)\nonumber \\
\dot{x}_{23} &=& {x_{23}}({x_{12}} {x_2} {x_{23}} ({k_2}^2+4 {k_2}{k_3}+{k_2}+2 {k_3}) (2 {x_{12}} {x_{13}}({k_2}+{k_3}-1)-{x_{12}}^2-{x_{13}}^2)+(1-{k_2}) {x_{12}}^2{x_{13}}\nonumber\\
&&  {x_2}^2({k_2}^2+3 {k_2} {k_3}+{k_2}+2{k_3})+{x_{12}} {x_2} {x_{23}}^3({k_2}^2+{k_2}+2{k_3})+{x_{23}}^2 (4 {x_{12}} {x_2} ({k_2}+{k_3}-1)({k_2} {x_{13}}\nonumber\\
&& +{k_3} {x_{12}})+({k_2}-2) ({k_2}-1) {k_2}{x_{12}}^2 {x_{13}}-({k_2}-1) {k_2} {x_{13}} {x_2}^2))
\end{eqnarray}  }
If we denote by $P_2(x_1, x_2, x_3, x_4)$ the $\dot{x}_2$, and by $P_{ij}(x_2, x_{12}, x_{13}, x_{23})$ the $\dot{x}_{ij}$ for $i< j=1,2,3$ respectively then the degree of vector field $X = P_2(x_2, x_{12}, x_{13}, x_{23})\frac{\partial}{\partial x_2}+ \sum_{i<j=1,2,3}P_{ij}(x_2, x_{12}, x_{13}, x_{23})\frac{\partial}{\partial x_{ij}}$ is $6$.  Now, in order
to study the singularities at infinity of (\ref{systemStiefel}), we should write this system in the local charts of the Poincar\'e compactification.  Because we are interested for positive values of $x_2, x_{12}, x_{13}$ and $x_{23}$ we study the behavior of the previous system only in to the chart $(U_1, \varphi_1)$.  Therefore from (\ref{Xartis11}) we have:
{\small \begin{eqnarray}\label{SystemStiefelP}
\dot{z}_1 &=& ({k_2}+1) {z_1} ({k_2}+2 {k_3}) (({z_1}-1) {z_2}^2 {z_3}(2 {z_1}-{k_2} ({z_1}-1))+{k_3} {z_1} ({z_1}^2 {z_2}-{z_1} {z_3}-{z_2} ({z_2}-{z_3})^2))\nonumber\\
\dot{z}_{2} &=&({k_2}+1) {z_2} ({k_2}+2 {k_3}) (-{z_3} ({z_1}^2(-2 {z_2} ({k_2}+{k_3}-1)+({k_2}-2){z_2}^2+{k_2}+{k_3}-1)+{z_2}^2)\nonumber\\
&& -{z_1} {z_2} {z_3}^2+{z_1} {z_2} ({z_2}-{z_1}) ({z_1}+{z_2}))\nonumber\\
\dot{z}_3 &=& (-{k_2}-1) {z_3} ({k_2}+2 {k_3}) ({k_2} {z_1} {z_2} ({z_1} {z_2}{z_3}+({z_1}-{z_2})^2-{z_3}^2)+{z_1}^2 ({k_3} {z_3}-2 {z_2}^2 ({k_3}+{z_3}-1))+{z_2}^2 {z_3})\nonumber\\
\dot{z}_4 &=& 2 {z_4} (-{z_2}^2 (({k_2}-2) {z_1}^2 {z_3} ({k_2}({k_2}+{k_3})+{k_3})+2 {z_1} ({k_2}+{k_3}-1) ({k_2}{z_3}+{k_3} {z_1})+{z_3} ({k_2}{k_3}+{k_2}+{k_3}))\nonumber\\
&& +{k_2} {k_3} {z_1} {z_2} (-2 {z_1} {z_3} ({k_2}+{k_3}-1)+{z_1}^2+{z_3}^2)-{k_3}{z_1}^2 {z_3} ({k_2} {k_3}+{k_2}+{k_3})+{k_2} {k_3}{z_1} {z_2}^3)
\end{eqnarray} }
\begin{theorem}
Let the Stiefel manifold $V_{1+k_2}\bb{R}^n \cong\SO(n)/\SO(k_3)$, with $n = 1+k_2+k_3$.   The normalized Ricci flow on the space of invariant Riemannian metrics on $\SO(1+2k_2)/\SO(k_2)$, for $k_2 \geq 6$ has at least four singularities at infinity.  
\end{theorem}
\begin{proof}
We have $k_2=k_3$ so the system (\ref{SystemStiefelP}) can be written as follows
\begin{eqnarray}\label{Stiefelk2}
&&\dot{z}_1  = 3 (-{k_2}-1) {k_2} {z_1} ({k_2} {z_1} ({z_1}^2(-{z_2})+{z_1} {z_3}+{z_2}({z_2}-{z_3})^2)+({z_1}-1) {z_2}^2 {z_3} ({k_2}({z_1}-1)-2 {z_1}))=0 \nonumber\\
&&\dot{z}_2 = 3 (-{k_2}-1) {k_2} {z_2} ({z_1}^2 {z_3} (({k_2}-2){z_2}^2+(2-4 {k_2}) {z_2}+2 {k_2}-1)+{z_1}^3{z_2}+{z_1} {z_2} ({z_3}^2-{z_2}^2)+{z_2}^2{z_3})=0 \nonumber\\
&&\dot{z}_3 = 3 (-{k_2}-1) {k_2} {z_3} ({k_2} {z_1} ({z_1}^2{z_2}+{z_1} ({z_2}^2 ({z_3}-4)+{z_3})+{z_2}^3-{z_2} {z_3}^2)+{z_2}^2({z_3}-2 {z_1}^2 ({z_3}-1)))=0.
\end{eqnarray}
We consider a polynomial ring $R = \bb{Q}[z, z_1, z_2, z_3]$ and a ideal $I$ generated by the polynomials $\{\dot{z}_1, \dot{z}_2, \dot{z}_3, z\, z_1\,$ $z_2\, z_3 -1\}$.  We take a lexicographic order $>$ with $z > z_3 > z_2 > z_1$ for a monomial ordering on $R$.  Then, by the aid of computer, we see that a Gr\"obner basis for the $I$ contains the polynomial $k_2(1+k_2)(z_1-1)F_(z_1)$ where $F_1(z_1)$ given by
{\footnotesize \begin{eqnarray*}
&&F_1(z_1) = 4 {k_2}^{12}+48 {k_2}^{11}-176 {k_2}^{10}+96 {k_2}^9+248 {k_2}^8-320 {k_2}^7+16 {k_2}^6+160 {k_2}^5-92 {k_2}^4+16 {k_2}^3+(289{k_2}^{12}\\
&& -2414 {k_2}^{11}+7267 {k_2}^{10}-7656 {k_2}^9-4925 {k_2}^8+16858 {k_2}^7-10115 {k_2}^6-1924 {k_2}^5+4556 {k_2}^4-4032{k_2}^3\\
&& +2864 {k_2}^2-832 {k_2}+64) {z_1}^{10}+(-3094{k_2}^{12}+21740 {k_2}^{11}-53420 {k_2}^{10}+43474 {k_2}^9+25150{k_2}^8-59760 {k_2}^7\\
&& +32856 {k_2}^6-11070 {k_2}^5+1676{k_2}^4+4944 {k_2}^3-3104 {k_2}^2+672 {k_2}-64){z_1}^9+(11951 {k_2}^{12}-62636 {k_2}^{11}+87727 {k_2}^{10}\\
&& +63904{k_2}^9-306211 {k_2}^8+405020 {k_2}^7-357771 {k_2}^6+244112{k_2}^5-115968 {k_2}^4+37280 {k_2}^3-9200 {k_2}^2+1664{k_2}\\
&& -128) {z_1}^8+(-23556 {k_2}^{12}+81980 {k_2}^{11}+1572{k_2}^{10}-381432 {k_2}^9+817168 {k_2}^8-1006372 {k_2}^7+869836{k_2}^6\\
&& -522872 {k_2}^5+212404 {k_2}^4-61192 {k_2}^3+12816{k_2}^2-1760 {k_2}+128) {z_1}^7+(27243 {k_2}^{12}-54638{k_2}^{11}-133224 {k_2}^{10}\\
&& +635796 {k_2}^9-1177180 {k_2}^8+1342266{k_2}^7-995716 {k_2}^6+467944 {k_2}^5-133163 {k_2}^4+19352{k_2}^3+376 {k_2}^2-672 {k_2}\\
&& +80) {z_1}^6+(-20366{k_2}^{12}+17416 {k_2}^{11}+186738 {k_2}^{10}-653354 {k_2}^9+1066178{k_2}^8-1003674 {k_2}^7+537954 {k_2}^6\\
&& -128838 {k_2}^5-18604{k_2}^4+23530 {k_2}^3-7564 {k_2}^2+1176 {k_2}-80){z_1}^5+(10457 {k_2}^{12}+3504 {k_2}^{11}-152744{k_2}^{10}\\
&& +427896 {k_2}^9-550638 {k_2}^8+347496 {k_2}^7-51536{k_2}^6-73272 {k_2}^5+54525 {k_2}^4-17048 {k_2}^3+2464 {k_2}^2-64{k_2}\\
&& -16) {z_1}^4+(-3592 {k_2}^{12}-7144 {k_2}^{11}+75380{k_2}^{10}-160756 {k_2}^9+138300 {k_2}^8-14956 {k_2}^7-63900{k_2}^6+50796 {k_2}^5\\
&& -15124 {k_2}^4+196 {k_2}^3+1016 {k_2}^2-232{k_2}+16) {z_1}^3+(768 {k_2}^{12}+3264 {k_2}^{11}-20788{k_2}^{10}+32264 {k_2}^9-10448 {k_2}^8\\
&& -20736 {k_2}^7+23488{k_2}^6-7528 {k_2}^5-1728 {k_2}^4+1904 {k_2}^3-508 {k_2}^2+48{k_2}) {z_1}^2+(-88 {k_2}^{12}-656 {k_2}^{11}+3024{k_2}^{10}\\
&& -3152 {k_2}^9-1552 {k_2}^8+4816 {k_2}^7-2496 {k_2}^6-624 {k_2}^5+1064 {k_2}^4-384 {k_2}^3+48 {k_2}^2) {z_1}.
\end{eqnarray*} }
If $z_1 = 1$, then from $\dot{z}_1=0$ in system (\ref{Stiefelk2})   we take $z_2 = z_3$ and $z_3 = (z_2^2-1)/z_2$ (this does not give a positive solution).  By substituting $z_2=z_3$ back into (\ref{Stiefelk2}) we find the real solutions 
$$
z_1 = 1, \  z_2 = z_3 =  \frac{-\sqrt{2 {k_2}^2-2 {k_2}+1}+2 {k_2}-1}{{k_2}-1} \ \mbox{and}\ z_1 = 1, \ z_2 = z_3 = \frac{\sqrt{2 {k_2}^2-2 {k_2}+1}+2 {k_2}-1}{{k_2}-1}.
$$ 
These are known as Jensen's Einstein metrics on Stiefel manifolds.
Now if $z_1\neq 1$ then $F_1(z_1) = 0$ and we claim that the equation $F_1(z_1) = 0$ has at least two positive roots.  We see that $F_1(1) = 16 {k_2}^{12}+464 {k_2}^{11}+1356 {k_2}^{10}-2920 {k_2}^9-3910{k_2}^8+10638 {k_2}^7-17384 {k_2}^6+16884{k_2}^5-10454{k_2}^4+4566 {k_2}^3-792 {k_2}^2 > 0$ and for $k_2 \geq 6$ we have 
$F_1(3) < 0$ and $F_1(6) > 0$.  Thus for $k_2\geq 6$ we have at least two solutions $z_1 = \alpha_1, \beta_1$ of $F_1(z_1) = 0$ with $1< \alpha_1 < 3$ and $3< \alpha_2 < 6$.  Next, we consider the ideal $J$ generated by the $\{\dot{z}_1, \dot{z}_2, \dot{z}_3, z\, z_1\,$ $z_2\, z_3\, (z_1-1) -1\}$ and take a lexicographic order $>$ with $z > z_2 > z_3 > z_1$ for a monomial ordering on $R$.  Then, by the aid of computer, we see that a Gr\"obner basis for the ideal $J$ contains the polynomials $F_1(z_1)$ and $q_3(k_2)z_3 = \sum_{\ell = 0}^{9}p_{\ell}(k_2)z_1^{\ell}$.  Also, for the same ideal $J$ and the lexicographic order $>$ with $z > z_3 > z_2 > z_1$ for monomials on $R$, we see that a Gr\"obner basis for $J$ contains the polynomial $q_2(k_2)z_2 = \sum_{\ell = 0}^{9}p_{\ell}(k_2)z_1^{\ell}$, where $q_i(k_2)$ $(i=2, 3)$ and $p_{\ell}(k_2)$ $(\ell = 0, 1,\ldots, 9)$ are polynomials of $k_2$ of degree $58$ and $60$ respectively with integer coefficients.  It is easy to see that $q_i(k_2) \neq 0$ $(i=2,3)$ for $k_2>1$.  Thus, for the positive solutions $z_1=\alpha_1, \beta_1$ found above we obtain real values $z_2 = \alpha_2, \beta_2$ and $z_3 = \al_3, \beta_3$ as solutions of the system (\ref{Stiefelk2}).  We  claim that $\al_2, \beta_2, \al_3, \beta_3$ are positive.  Take a lexicographic order $>$ with $z > z_1 > z_3 > z_2$ for a monomial ordering on the ring $R$. Then, a Gr\"obner basis for the ideal $J$ contains a polynomial $F_2(z_2)$ of $z_2$, which is given by
{\footnotesize \begin{eqnarray*}
&& F_2(z_2) = (2680947582 ({k_2}-6)^7+14591972187 ({k_2}-6)^6+59640424960 ({k_2}-6)^5+180449511648 ({k_2}-6)^4\\
&& +392350767872 ({k_2}-6)^3+580027773696 ({k_2}-6)^2+522221572096 ({k_2}-6)+216162222080)z_2^{10}\\
&&+(-50225598762 ({k_2}-6)^7-276982406346 ({k_2}-6)^6-1149224425652({k_2}-6)^5-3535515192656 ({k_2}-6)^4\\
&& -7827543360064({k_2}-6)^3-11798003601664 ({k_2}-6)^2-10842357528576({k_2}-6)-4585758048256)z_2^9\\
&&+(390159654498 ({k_2}-6)^7+2172262119880 ({k_2}-6)^6+9119121333384({k_2}-6)^5+28438722948740 ({k_2}-6)^4\\
&& +63932956074624 ({k_2}-6)^3+97997952225152 ({k_2}-6)^2+91718708609024({k_2}-6)+39559071548416)z_2^8\\
&&+(-1636143378304 ({k_2}-6)^7-9149184702208 ({k_2}-6)^6-38667393707288({k_2}-6)^5-121655123121672 ({k_2}-6)^4\\
&& -276434836684176({k_2}-6)^3-429032167222784 ({k_2}-6)^2-407238153983232({k_2}-6)-178416492097536)z_2^7\\
&&+(4077237826688 ({k_2}-6)^7+22741035201568 ({k_2}-6)^6+96082339120904({k_2}-6)^5+302818376659728 ({k_2}-6)^4\\
&& +690563967560064({k_2}-6)^3+1077488181695360 ({k_2}-6)^2+1029907447940096({k_2}-6)+455097671436288)z_2^6\\
&&+(-6367615594352 ({k_2}-6)^7-35174632043952 ({k_2}-6)^6-147420550162112({k_2}-6)^5-461548083948832 ({k_2}-6)^4\\
&& -1046969395743360({k_2}-6)^3-1626943248550784 ({k_2}-6)^2-1550563723207680({k_2}-6)-683916568270848)z_2^5\\
&&+(6500710198416 ({k_2}-6)^7+35468514250320 ({k_2}-6)^6+146902006546624({k_2}-6)^5+454750667573776 ({k_2}-6)^4\\
&& +1020485202419072({k_2}-6)^3+1569594503117888 ({k_2}-6)^2+1481383976485632({k_2}-6)+647378577709056)z_2^4\\
&&+(-4376025268608 ({k_2}-6)^7-23628587384704 ({k_2}-6)^6-96812750676960 ({k_2}-6)^5-296392431072480 ({k_2}-6)^4\\
&& -657659598211392({k_2}-6)^3-1000042532332416 ({k_2}-6)^2-933014709911808({k_2}-6)-403027461872640)z_2^3\\
&&+(1873316885024 ({k_2}-6)^7+10053816339984 ({k_2}-6)^6+40911319513376({k_2}-6)^5+124311949165760 ({k_2}-6)^4\\
&& +273618707342592({k_2}-6)^3+412536333595392 ({k_2}-6)^2+381469386539520({k_2}-6)+163261397836800)z_2^2\\
&&+(-460264120480 ({k_2}-6)^7-2465535026848 ({k_2}-6)^6-10005481382592({k_2}-6)^5-30298379604480 ({k_2}-6)^4\\
&& -66422377562112({k_2}-6)^3-99697425799680 ({k_2}-6)^2-91739835878400({k_2}-6)-39057790464000)z_2\\
&& 49150850208 ({k_2}-6)^7+263520189632 ({k_2}-6)^6+1069436792448({k_2}-6)^5+3236342259456 ({k_2}-6)^4\\
&& +7086348357120({k_2}-6)^3+10618426598400 ({k_2}-6)^2+9750544128000({k_2}-6)+4141186560000
\end{eqnarray*} }
Then we see that, for $k_2\geq 6$, the coefficients of the polynomial $F_2(z_2)$ are positive for even degree and negative for odd degree terms.  Thus if the equation $F_2(z_2)= 0 $ has real solutions,  then these are all positive.  In particular $z_2 = \al_2, \beta_2$ are positive.  Now we take a lexicographic order $>$ with $z >z_2 > z_1 > z_3$ for a monomial ordering on $R$.  Then, by the aid of computer, we see that a Gr\"obner basis for the ideal $J$ contains the polynomial$F_3(z_3$), where the polynomial $F_3(z_3)$ can be written as
{\footnotesize \begin{eqnarray*}
&& F_3(z_3) = (21286672726062828160 ({k_2}-6)^7+58254716280111893248({k_2}-6)^6+127455519890994757632 ({k_2}-6)^5\\
&&+217655806674022477824({k_2}-6)^4+279537041132943605760 ({k_2}-6)^3+253958157960834252800({k_2}-6)^2\\
&&+145499460899831808000 ({k_2}-6)+39531747173990400000)z_3^{10} +(-292777902042514008608 ({k_2}-6)^7\\
&&-814588260579033355648 ({k_2}-6)^6-1811976784418141538304 ({k_2}-6)^5-3146205600377563549696({k_2}-6)^4\\
&& -4108991272662490095616 ({k_2}-6)^3-3796725395351760764928({k_2}-6)^2-2212811038797252526080 ({k_2}-6)\\
&& -611729579268322099200 )z_3^9 + (1660257266722947751552 ({k_2}-6)^7+4695714013629690519424({k_2}-6)^6\\
&& +10616322476744264890880 ({k_2}-6)^5+18734391448720331455488 ({k_2}-6)^4+24867099690573745029120 ({k_2}-6)^3\\
&& +23354475414681440354304({k_2}-6)^2+13836689687674552975360 ({k_2}-6)+3889096495709532979200 )z_3^8\\
&& +(-5079593933368097151208 ({k_2}-6)^7-14594154103096371015040({k_2}-6)^6-33503861421156201380608 ({k_2}-6)^5\\
&& -60018579090397951259648({k_2}-6)^4-80858976481531529269248 ({k_2}-6)^3-77073028202750859501568({k_2}-6)^2\\
&& -46344650438235365801984 ({k_2}-6)-13221646074705577377792)z_3^7 + (9308466298136566608492 ({k_2}-6)^7\\
&& +27115170479328586830320({k_2}-6)^6+63066218770915166201344 ({k_2}-6)^5+114399455055489074824960({k_2}-6)^4\\
&& +156004137935460456838144 ({k_2}-6)^3+150475504061699408097280({k_2}-6)^2+91547705640188726935552 ({k_2}-6)\\
&& +26423046593985210351616 )z_3^6 + (-11095597642774291380236 ({k_2}-6)^7-32642462593282533158032({k_2}-6)^6
\end{eqnarray*} }
{\footnotesize \begin{eqnarray*}
&& -76608697711598428869856 ({k_2}-6)^5-140125082414737357727872({k_2}-6)^4-192577122303638985157632 ({k_2}-6)^3\\
&& -187123466429387154372608({k_2}-6)^2-114646689401743117639680 ({k_2}-6)-33314932853414556467200 )z_3^5\\
&& +(9347621052101913097360 ({k_2}-6)^7+27635109722202090645504({k_2}-6)^6+65135190023242223199008 ({k_2}-6)^5\\
&& +119592275581440467432704({k_2}-6)^4+164920911373642400207616 ({k_2}-6)^3+160749870833420529548288({k_2}-6)^2\\
&& +98770860152536885063680 ({k_2}-6)+28778321592849865523200 )z_3^4 +(-5615419680567151075808 ({k_2}-6)^7\\
&& -16600290237902267697984({k_2}-6)^6-39112773412702956110336 ({k_2}-6)^5-71771758931254298414464({k_2}-6)^4\\
&& -98898714775619171369216 ({k_2}-6)^3-96307590510023343088640({k_2}-6)^2-59111874345673483827200 ({k_2}-6)\\
&& -17202705345352952832000 )z_3^3 +(2402583297061148772864 ({k_2}-6)^7+7056849843649244636416({k_2}-6)^6\\
&& +16525930989063218544768 ({k_2}-6)^5+30148489433708152475392({k_2}-6)^4+41309556136774635555840 ({k_2}-6)^3\\
&& +40006454668872611788800({k_2}-6)^2+24423031403576407296000 ({k_2}-6)+7069849819878827520000 )z_3^2\\
&& +(-688999017327743752704 ({k_2}-6)^7-1994617060599839321088({k_2}-6)^6-4609779004452472829952 ({k_2}-6)^5\\
&& -8307601122406729420800({k_2}-6)^4-11253654298468482048000 ({k_2}-6)^3-10781346008039639040000({k_2}-6)^2\\
&& -6514118885048832000000 ({k_2}-6)-1867029805159833600000 )z_3 + 109572943236302618624 ({k_2}-6)^7\\
&& +309741220607312068608({k_2}-6)^6+700554453138220400640 ({k_2}-6)^5+1237766738525319168000({k_2}-6)^4\\
&& +1646233648197304320000 ({k_2}-6)^3+1550332512446054400000({k_2}-6)^2+921701573594726400000 ({k_2}-6)\\
&&+260150759534592000000.
\end{eqnarray*} }
Then, we see that, for $k_2\geq 6$, the coefficients of the polynomial $F_3(z_3)$ are positive for even degree terms and negative for odd degree terms. Thus if the equation $F_3(z_3) = 0$ has real solutions, then these are all positive.  In particular, the solutions $z_3=\al_3, \beta_3$ are positive, thus we get our claim.  
\end{proof}

In order to find the singularities at infinity of the above system we set $z_4 = 0$.

Next, we compute the fixed points of the system ${\dot{z}_1 = 0, \dot{z}_2 = 0, \dot{z}_3 = 0}$ for specific values of $k_2$ and $k_3>1$.  We have:
 
 
\smallskip 
\noindent$\bullet$\ $V_{5}\bb{R}^7$:
$(4.1466, 4.07919, 1.03361)$,  $(2.29783, 3.43436, 3.98856)$,  $(1, 2.54858, 2.54858)$, $(1, 0.78475, 0.78475)$. 

\smallskip
\noindent$\bullet$\ $V_{5}\bb{R}^8$: 
 $(5.39567, 4.8672, 2.16024)$,  $(2.31234,  4.49843,  4.93295)$, $(1,  3.29099, 3.29099)$, $(1,  0.709006, 0.709006)$. 


\centerline{$\ldots$}


\smallskip
\noindent$\bullet$\  $V_6\bb{R}^8$:
$(3.19365, 3.15771,  0.674502)$,  $(1.86343, 2.64311,  3.07833)$,   $(1,  2.20711,  2.20711)$, $(1,  0.792893,  0.792893)$.  

\smallskip
\noindent$\bullet$\ $V_6\bb{R}^9$:
$(3.99996, 3.71213,  1.41708)$,  $(1.89382, 3.36866, 3.73723)$,  $(1,  2.78078, 2.78078)$ $(1,  0.719224,  0.719224)$.
 

\centerline{$\ldots$}


\smallskip
\noindent$\bullet$\ $V_7\bb{R}^9$:
$(2.71186, 2.68928, 0.499721)$,  $(1.64442, 2.25706, 2.6166)$, $(1,  2, 2)$, $(0.95544,$ $0.798009,$ $0.734193)$,\\   $(0.805105,$ $0.771014,$ $0.379868)$,  $(1,  0.8, 0.8)$. 

\smallskip
\noindent$\bullet$\ $V_7\bb{R}^{10}$:
 $(3.30651, 3.12526, 1.05079)$,  $(1.67763, 2.81537,  3.13489)$, $(1,  2.47178,  2.47178)$, $(1,  0.72822,  0.72822)$.
 

\centerline{$\ldots$}

 
\smallskip 
\noindent$\bullet$\ $V_8\bb{R}^{10}$:
$(2.41937, 2.40377, 0.396819)$, $(1.51286, 2.02874, 2.33539)$, $(1, 1.86038, 1.86038)$, $(1, 0.806287, 0.806287)$,\\ $(0.98506, 0.805755, 0.78523)$, $(0.791817, 0.770023, 0.312754)$.

\smallskip
\noindent$\bullet$\ $V_8\bb{R}^{11}$: 
$(2.89171, 2.76727, 0.833614)$,  $(1.54531, 2.48719, 2.76871)$,  $(1, 2.26376, 2.26376)$, $(1, 0.736237, 0.736237)$.


\centerline{$\ldots$}






  




\noindent
From the above results we have: 

\begin{theorem}
The normalized Ricci flow on the space of invariant Riemannian metrics on $V_5\bb{R}^7$, possesses exactly four singularities at infinity.  These fixed points correspond (up to scale) to the $G$-invariant Einstein metrics on $V_5\bb{R}^7$ (cf. Theorem 4.6). 
\end{theorem}

We can also make the following:

\begin{conjecture}
Let $G/H$ be the Stiefel manifold $V_{1+k_2}\bb{R}^{n}$, with $n = 1+k_2+k_3$.  Then, for $k_2 \geq 4$ and $k_3 > 1$, $G/H$ has four singularities at infinity, and for $k_2\geq 6$ and $k_3 = 2$  it has six singularities at infinity.  These fixed points correspond to the $G$-invariant Einstein metrics on $G/H$.  
\end{conjecture}

\subsection{Ricci flow for the generalized flag manifolds} 
We study the behavior of the normalized Ricci flow for the generalized flag manifold $G/K$ with four isotropy summands and $b_2(G/K)=1$.  For this case we take the Ricci components (\ref{compI}) and scalar curvature (\ref{scalFlag}) of the metric (\ref{metricFlag}).  Then, system (\ref{normalRicci}) reduces to 

\begin{eqnarray*}
&&\dot{x}_1 = - \frac{1}{2 {d_1} {x_1} {x_2}^2 {x_3} {x_4}N}  \Big\lbrace {x_2} {x_4} ({x_3} ({A_{112}} {x_2}^2 (3 {d_1}+2({d_2}+{d_3}+{d_4}))+2{d_1} {x_1}^2 ({A_{112}}-{d_2})\\
&& -2{d_1} {x_1} {x_2} (2 {d_1}+{d_2}+{d_3}+{d_4}))+2{A_{123}} {x_1} (2 {d_1}({x_2}^2+{x_3}^2)+(-{d_2}-{d_3}-{d_4})({x_1}^2-{x_2}^2-{x_3}^2))\\
&& -2 {d_1} {d_3}{x_1}^2 {x_2})+2 {x_1} {x_2}^2 ({A_{134}} {x_3}^2 (2{d_1}+{d_2}+{d_3}+{d_4})-{A_{134}} {x1}^2
({d_2}+{d_3}+{d_4})\\
&&+{d_1} {x_1} {x_3}({A_{224}}-{d_4}))+{x_1} {x_4}^2 (2 {A_{134}} {x_2}^2(2 {d_1}+{d_2}+{d_3}+{d_4})+{A_{224}} {d_1} {x_1} {x_3}) \Big\rbrace             
\\
&&\dot{x}_2 = \frac{1}{2 {d_2} {x_1}^2 {x_2} {x_3} {x_4}N} \Big\lbrace {x_2} {x_4} ({x_3} (-2 {x_1}^2 ({A_{112}}-{d_2})({d_1} +2 {d_2}+{d_3}+{d_4})+{A_{112}} {x_2}^2({d_1} +{d_3}+{d_4})\\
&&+2 {d_1} {d_2} {x_1} {x_2})+2{x_1} (-{A_{123}} {x_1}^2 ({d_1}+2{d_2}+{d_3}+{d_4})+{A_{123}} {x_2}^2({d_1}+{d_3}+{d_4})+{d_2} {d_3} {x_1} {x_2})\\
&& -2{A_{123}} {x_1} {x_3}^2 ({d_1}+2 {d_2}+{d_3}+{d_4}))-2{A_{134}} {d_2} {x_1} {x_2}^2({x_1}^2+{x_3}^2+{x_4}^2)\\
&& -{A_{224}} {x_1}^2 {x_3}{x_4}^2 (2 {d_1}+3 {d_2} +2 ({d_3}+{d_4}))+2 {d_2} {x_1}^2{x_2}^2 {x_3} ({d_4}-{A_{224}}) \Big\rbrace 
\\
&&\dot{x}_3 = -\frac{1}{2 {d_3} {x_1}^2 {x_2}^2 {x_4} N} \Big\lbrace {x_2}{x_4} ({d_3} {x_3} ({A_{112}} (2{x_1}^2+{x_2}^2)-2 {x_1} ({d_1} {x_2}+{d_2}{x_1}))+2 {x_1} ({A_{123}}({x_1}^2 +{x_2}^2)\\
&&-{d_3}{x_1} {x_2})({d_1}+{d_2}+2 {d_3}+{d_4})-2 {A_{123}} {x_1} {x_3}^2({d_1}+{d_2}+{d_4}))+2 {x_1} {x_2}^2 ({A_{134}}{x_1}^2 ({d_1}+{d_2} \\
&& +2 {d_3}+{d_4})-{A_{134}} {x_3}^2({d_1}+{d_2}+{d_4})+{d_3} {x_1} {x_3}({A_{224}}-{d_4}))+{x_1} {x_4}^2 (2 {A_{134}} {x_2}^2({d_1}+{d_2}\\
&& +2 {d_3}+{d_4}) +{A_{224}} {d_3} {x_1}{x_3}) \Big\rbrace 
\\
&&\dot{x}_4 = \frac{1}{2 {d_4} {x_1}^2 {x_2}^2 {x_3} N} \Big\lbrace -{d_4} {x_2} {x_4} ({x_3} (2 {x_1}^2({A_{112}}-{d_2})+{A_{112}} {x_2}^2-2{d_1} {x_1}{x_2}) +2 {A_{123}} {x_1}({x_1}^2+{x_2}^2\\
&& +{x_3}^2)-2 {d_3} {x_1}^2{x_2})-2 {x_1} {x_2}^2 ({d_1}+{d_2}+{d_3}+2 {d_4})({A_{134}} ({x_1}^2+{x_3}^2)+{x_1} {x_3}({A_{224}}-{d_4}))\\
&& +{x_1} {x_4}^2 ({d_1}+{d_2}+{d_3})(2 {A_{134}} {x_2}^2+{A_{224}} {x_1} {x_3}) \Big\rbrace ,
\end{eqnarray*}
where $N = d_1+d_2+d_3+d_4$.  The above system is not polynomial, hence in order to apply the Poincar\'e compactification, we multiply it by the factor $x_1^2x_2^2x_3x_4d_1d_2d_3d_4N$.  After that we take the following system:
{ \begin{eqnarray}\label{systemFlag}
\dot{x}_1 &=& -{d_2} {d_3} {d_4} {x_1} ({x_2} {x_4}({x_3}({A_{112}} {x_2}^2 (3 {d_1}+2 ({d_2}+{d_3}+{d4}))+2
{d_1} {x_1}^2 ({A_{112}}-{d_2})-2{d_1} {x_1} {x_2} (2{d_1}+{d_2}\nonumber\\
&& +{d_3}+{d_4}))+2 {A_{123}} {x_1} (2{d_1} ({x_2}^2+{x_3}^2)+(-{d_2}-{d_3}-{d_4})({x_1}^2-{x_2}^2-{x_3}^2))-2 {d_1} {d_3}{x_1}^2 {x_2})\nonumber\\
&& +2 {x_1} {x_2}^2 ({A_{134}} {x_3}^2 (2{d_1}+{d_2}+{d_3}+{d_4})-{A_{134}} {x_1}^2({d_2}+{d_3}+{d_4})+{d_1} {x_1}{x_3}({A_{224}}-{d_4}))\nonumber\\
&& +{x_1} {x_4}^2 (2 {A_{134}} {x_2}^2(2 {d_1}+{d_2}+{d_3}+{d_4})+{A_{224}} {d_1}{x_1}{x_3}))
\nonumber\\
\dot{x}_2 &=& {d_1} {d_3} {d_4} {x_2} ({x_2} {x_4} ({x_3}(-2 {x_1}^2 ({A_{112}}-{d_2}) ({d_1}+2{d_2}+{d_3}+{d_4})+{A_{112}} {x_2}^2({d_1}+{d_3}+{d_4})\nonumber\\
&& +2 {d_1} {d_2}{x_1}{x_2}) +2{x_1} (-{A_{123}} {x_1}^2 ({d_1}+2{d_2}+{d_3}+{d_4})+{A_{123}} {x_2}^2({d_1}+{d_3}+{d_4})+{d_2} {d_3} {x_1} {x_2})\nonumber\\
&&-2{A_{123}} {x_1} {x_3}^2 ({d_1} +2{d_2} +{d_3}+{d_4}))-2 {A_{134}} {d_2}{x_1}{x_2}^2({x_1}^2+{x_3}^2+{x_4}^2)-{A_{224}} {x_1}^2 {x_3}{x_4}^2 (2 {d_1}\nonumber\\
&&+3 {d_2}+2 ({d_3}+{d_4})) +2 {d_2}{x_1}^2{x_2}^2 {x_3} ({d_4}-{A_{224}}))
\nonumber\\
\dot{x}_3 &=& -{d_1} {d_2} {d_4} {x_3} ({x_2} {x_4} ({d_3}{x_3} ({A_{112}} (2 {x_1}^2+{x_2}^2)-2 {x_1}({d_1} {x_2}+{d_2} {x_1}))+2 {x_1} ({A_{123}}({x_1}^2+{x_2}^2)\nonumber\\
&& -{d_3} {x_1} {x_2}) ({d_1} +{d_2}+2 {d_3}+{d_4})-2 {A_{123}} {x_1} {x_3}^2({d_1}+{d_2}+{d_4}))+2 {x_1} {x_2}^2 ({A_{134}}{x_1}^2 ({d_1}+{d_2}\nonumber\\
&& +2 {d_3}+{d_4})-{A_{134}} {x_3}^2({d_1} +{d_2}+{d_4})+{d_3} {x_1} {x_3}({A_{224}}-{d_4}))+{x_1} {x_4}^2 (2 {A_{134}} {x_2}^2({d_1}+{d_2}\nonumber\\
&& +2 {d_3}+{d_4})+{A_{224}} {d_3} {x_1}{x_3}))\nonumber\\
\dot{x}_4 &=& {d_1} {d_2} {d_3} {x_4} (-{d_4} {x_2} {x_4}({x_3} (2 {x_1}^2 ({A_{112}}-{d_2})+{A_{112}}{x_2}^2-2{d_1} {x_1} {x_2})+2 {A_{123}} {x_1}({x_1}^2+{x_2}^2+{x_3}^2)\nonumber\\
&& -2 {d_3} {x_1}^2{x_2})-2 {x_1} {x_2}^2 ({d_1}+{d_2}+{d_3}+2 {d_4})({A_{134}} ({x_1}^2+{x_3}^2)+{x_1} {x_3}({A_{224}}-{d_4}))\nonumber\\
&& +{x_1} {x_4}^2 ({d_1}+{d_2}+{d_3})(2 {A_{134}} {x_2}^2+{A_{224}} {x_1} {x_3})).
\end{eqnarray} }

If we denote by $P_i(x_1, x_2, x_3, x_4)$ the $\dot{x}_i$, for $i=1,2,3,4$ respectively then the degree of vector field $X = \sum_{i=1}^4P_i(x_1, x_2, x_3, x_4)\frac{\partial}{\partial x_i}$ is $6$.  Next,
it is easy to see the following lemma
\begin{lemma}
The coordinate planes along with the straight line $\gamma(t) = (t, 2t, 3t, 4t)$, remain invariant under the normalized Ricci flow defined by the system (\ref{systemFlag}).
\end{lemma}

In order to study the singularities at infinity of (\ref{systemFlag}), we should write this system in the local charts of the Poincar\'e compactification.  Since we are interested for positive values of $x_1, x_{2}, x_{3}$ and $x_{4}$, we study the behavior of the previous system only in to the chart $(U_1, \varphi_1)$.  Therefore from (\ref{Xartis11}) we have:
{\small \begin{eqnarray}\label{PoincareFlag}
\dot{z}_1 &=& {d_3} {d_4} {z_1} ({d_1}+{d_2}+{d_3}+{d_4}) ({z_3}({A_{112}} {z_1} {z_2} ({z_1}^2 ({d_1}+2 {d_2})-2{d_1})+2 {A_{123}} {z_1} (({z_1}^2-1)({d_1}+{d_2})\nonumber\\
&& +{z_2}^2 ({d_2}-{d_1}))-2{d_1} {z_2}({A_{224}} {z_3}+{d_2} ({z_1}-1) {z_1}))+2 {A_{134}}{d_2} {z_1}^2 ({z_2}^2+{z_3}^2-1))
\nonumber\\
\dot{z}_2 &=& 2 {d_2} {d_4} {z_1} {z_2} ({d_1}+{d_2}+{d_3}+{d_4})({z_3} ({d_3} {z_1} ({A_{112}} {z_1}{z_2}+{d_1}
(-{z_2})+{d_1})-{A_{123}} ({z_1}^2({d_1}-{d_3})+{d_1}+{d_3}) \nonumber\\
&& +{A_{123}} {z_2}^2({d_1}+{d_3}))+{A_{134}} {z_1} (({z_2}^2-1)({d_1}+{d_3})+{z_3}^2 ({d_3}-{d_1}))) 
\nonumber\\
\dot{z}_3 &=& -{d_2} {d_3} {z_3} ({d_1}+{d_2}+{d_3}+{d_4}) (-2{d_4} {z_1} {z_3} ({z_1} {z_2} ({A_{112}}{z_1}-{d_1})+{A_{123}} ({z_1}^2+{z_2}^2-1))  \nonumber\\
&& +2{z_1}^2 ({A_{134}} ({z_2}^2({d_1}-{d_4})+{d_1}+{d_4})+{d_1} {z_2}
({A_{224}}-{d_4}))-{z_3}^2 (2 {A_{134}} {z_1}^2({d_1}+{d_4})+{A_{224}} {d_1} {z_2}))\nonumber\\
\dot{z}_4 &=& {d_2} {d_3} {d_4} {z_4} ({z_1} {z_3} ({z_2}({A_{112}} {z_1}^2 (3 {d_1}+2 ({d_2}+{d_3}+{d_4}))+2{A_{112}} {d_1}-2 {d_1} {z_1} (2{d_1}+{d_2}+{d_3}+{d_4})-2 {d_1} {d_2})\nonumber\\
&& +2 {A_{123}}{z_1}^2 (2 {d_1}+{d_2}+{d_3}+{d_4})+2 {A_{123}} {z_2}^2 (2{d_1}+{d_2}+{d_3}+{d_4})-2 {A_{123}}({d_2}+{d_3}+{d_4})-2 {d_1} {d_3} {z_1})\nonumber\\
&& +2 {z_1}^2({A_{134}} {z_2}^2 (2 {d_1}+{d_2}+{d_3}+{d_4})-{A_{134}}({d_2}+{d_3}+{d_4})+{d_1} {z_2}({A_{224}}-{d_4}))\nonumber\\
&& +{z_3}^2 (2 {A_{134}} {z_1}^2 (2{d_1}+{d_2}+{d_3}+{d_4})+{A_{224}} {d_1} {z_2})).
\end{eqnarray} }

We set $z_4 = 0$ in order to obtain the behavior at infinity of the system (\ref{PoincareFlag}).  Next
we will study the system $\{\dot{z}_1 = 0, \dot{z}_2 = 0, \dot{z}_3 = 0\}$, in any case of flag manifolds with $b_2(G/K)$ separately.  First we substitute the values of the dimensions $d_i$, $i=1,2,3,4$ and the numbers $A_{224}, A_{112}, A_{123}$ and $A_{134}$ from Table 4.  It is easy to see that system (\ref{PoincareFlag}) has always a singularity located at $(2,3,4)$, which corresponds to the K\"ahler metric $(1,2,3,4)$.  For the flag manifolds which correspond to the exceptional Lie groups $\F_4, \E_7$ and $\E_8(\al_3)$ we found two more fixed points and for $\E_{8}(\al_6)$ four more. Actually we have:
    
\noindent $\blacktriangleright$\ $\F_4/\SU(3)\times\SU(2)\times\U(1)$
\begin{eqnarray*}
\dot{z}_1 &=& 11520 {z_1} ({z_1}^3 (8 {z_2}+5) {z_3}+2 {z_1}^2({z_2}^2-18 {z_2} {z_3}+{z_3}^2-1)+{z_1} ({z_2}
({z_2}+32)-5) {z_3}-4 {z_2} {z_3}^2) \\
\dot{z}_2 &=& 23040 {z_1} {z_2}(3 ({z_2}-1) {z_3} ({z_1}^2-6{z_1}+2 {z_2}+2)+4 {z_1} ({z_2}^2-1)-2 {z_1}{z_3}^2)\\
\dot{z}_3 &=& -2880 {z_3} (-12 {z_1}{z_3} ({z_1}^2+2 ({z1}-6){z_1} {z_2}+{z_2}^2-1)-24 {z_3}^2({z_1}^2+{z_2})+8 {z_1}^2 (({z_2}-12) {z_2}+3))
\end{eqnarray*}
The solutions are: $(0.970488, 0.229171, 1.0097)$ and $(1.27614, 1.95786, 2.31788)$

\noindent $\blacktriangleright$\ $\E_7/\SU(4)\times\SU(3)\times\SU(2)\times\U(1)$
\begin{eqnarray*}
\dot{z}_1 &=& 976896 {z_1} ({z_1}^3 (10 {z_2}+7) {z_3}+{z_1}^2({z_2}^2-36 {z_2} {z_3}+{z_3}^2-1)-{z_1}(({z_2}-28) {z_2}+7) {z_3}-2 {z_2} {z_3}^2)\\
\dot{z}_2 &=& 1953792 {z_1} {z_2} (3 ({z_2}-1) {z_3} ({z_1}^2-6{z_1}+2 {z_2}+2)+2 {z_1} ({z_2}^2-1)-{z_1}{z_3}^2)\\
\dot{z}_3 &=& 976896 {z_3} (3 {z_1} {z_3} ({z_1}^2+2 ({z_1}-6) {z_1}{z_2}+{z_2}^2-1)+3{z_3}^2 (3 {z_1}^2+2{z_2})+{z_1}^2 (-({z_2}-3)) (7 {z_2}-3))
\end{eqnarray*}
The solutions are: $(0.823351, 1.29423, 1.34989)$ and $(0.991279, 0.578307, 1.13127)$

\noindent $\blacktriangleright$\ $\E_8(\al_6)/\SU(7)\times\SU(2)\times\U(1)$
\begin{eqnarray*}
&&\dot{z}_1 = 15059072 {z_1} ({z_1}^3 (16 {z_2}+11) {z_3}+2 {z_1}^2({z_2}^2-30 {z_2} {z_3}+{z_3}^2-1)-{z_1}(({z_2}-48) {z_2}+11) {z_3}-4 {z_2} {z_3}^2)\\
&& \dot{z}_2 = 30118144 {z_1} {z_2} (5 ({z_2}-1) {z_3} ({z_1}^2-6{z_1}+2 {z_2}+2)+4 {z_1} ({z_2}^2-1)-2 {z_1}{z_3}^2)\\
&& \dot{z}_3 = 15059072 {z_3} (5 {z_1} {z_3} ({z_1}^2+2 ({z_1}-6){z_1} {z_2}+{z_2}^2-1)+2 {z_3}^2 (7 {z_1}^2+5{z_2})-2 {z_1}^2 (5 ({z_2}-4) {z_2}+7))
\end{eqnarray*}
The solutions are: $(0.91333, 1.41368, 1.51968)$ and $(0.966311, 0.489832, 1.08091)$

\noindent $\blacktriangleright$\ $\E_8(\al_3)/\SO(10)\times\SU(3)\times\U(1)$
\begin{eqnarray*}
&&\dot{z}_1 =7151616 {z_1} ({z_1}^3 (18 {z_2}+13) {z_3}+{z_1}^2({z_2}^2-60 {z_2} {z_3}+{z_3}^2-1)+{z_1} ((44-3
{z_2}) {z_2}-13) {z_3}-2 {z_2} {z_3}^2)\\
&& \dot{z}_2 = 14303232 {z_1} {z_2} (5 ({z_2}-1) {z_3} ({z_1}^2-6{z_1}+2 {z_2}+2)+2 {z_1} ({z_2}^2-1)-{z_1}{z_3}^2)\\
&& \dot{z}_3 = 7151616 {z_3}(5 {z_1} {z_3} ({z_1}^2+2 ({z_1}-6){z_1} {z_2}+{z_2}^2-1)+{z_3}^2 (17 {z_1}^2+10{z_2})+{z_1}^2 (5 (8-3 {z_2}) {z_2}-17))
\end{eqnarray*}
The solutions are: $(0.649612, 1.10943, 1.06103)$, $(0.763357, 1.00902, 0.191009)$, $(1.15607, 1.01783, 0.214618)$ and $(1.09705, 0.770347, 1.29696)$

We work on the chart $(U_1, \varphi_1)$ that corresponds to the plane $y_1 = 1$, so we consider metrics whose are the form $(1, \alpha, \beta, \gamma)$, where $\alpha, \beta$ and $\gamma$ are the solution of the system (\ref{PoincareFlag}) for $z_4=0$.  These metrics are invariant Einstein metrics, and the one with coefficients $(1, 2, 3, 4)$ is the unique K\"ahler-Einstein that admits $M$.  We have thus proved the following.

\begin{theorem}
Let $M = G/K$ be a generalized flag manifold with four isotropy summands and $b_2(M) = 1$.  The normalized Ricci flow, on the space of invariant Riemannian metrics on $M$, possesses exactly three singularities at infinity in case of $\F_4$, $\E_7$, $\E_8(\al_6)$ and exactly five in case of $\E_8(\al_3)$.  
These fixed points correspond to the $G$-invariant Einstein metrics on $M$ (cf. Theorem \ref{TheoremFlag}).
\end{theorem}

\medskip
\noindent 
{\bf Acknowledgement. }This work was developed and completed while the author had a research grant from DAAD  at Philipps-Universit\"at Marburg, during Fall 2019.

\end{document}